\documentclass[reqno,10pt]{amsart}

\usepackage{times}
\usepackage[T1]{fontenc}
\usepackage{mathrsfs}

\usepackage{latexsym}
\usepackage[dvips]{graphics}
\usepackage{epsfig}
\usepackage{amsmath,amsfonts,amsthm,amssymb,amscd}
\input amssym.def
\input amssym.tex
\usepackage{color}
\usepackage{graphicx}
\usepackage{verbatim}

\renewcommand{\i}{i} 

\newcommand{\mattwo}[4]
{\left(\begin{array}{cc}
                        #1  & #2   \\
                        #3 &  #4
                          \end{array}\right) }

\newcommand{\glfour}{{\rm GL}(4)}
\newcommand{\glsix}{{\rm GL}(6)}

\newcommand{\soe}{{\rm SO(even)}}
\newcommand{\soo}{{\rm SO(odd)}}

\newcommand{\sy}{{\rm USp}}
\newcommand{\un}{{\rm U}}
\newcommand{\so}{{\rm O}} 

\newcommand\be{\begin{equation}}
\newcommand\ee{\end{equation}}
\newcommand\bea{\begin{eqnarray}}
\newcommand\eea{\end{eqnarray}}
\newcommand\bi{\begin{itemize}}
\newcommand\ei{\end{itemize}}
\newcommand\ben{\begin{enumerate}}
\newcommand\een{\end{enumerate}}
\newcommand\bc{\begin{center}}
\newcommand\ec{\end{center}}
\newcommand\ba{\begin{array}}
\newcommand\ea{\end{array}}

\newcommand{\R}{\ensuremath{\mathbb{R}}}

\newcommand{\Z}{\ensuremath{\mathbb{Z}}}
\newcommand{\Q}{\mathbb{Q}}
\newcommand{\T}{\mathbb{T}}

\newtheorem{thm}{Theorem}[section]
\newtheorem{conj}[thm]{Conjecture}

\theoremstyle{definition}
\newtheorem{rek}[thm]{Remark}

\newcommand{\gl}{\lambda}

\newcommand{\hphi}{\widehat{\phi}}  

\newcommand{\foh}{\frac{1}{2}}  

\newcommand{\plogce}{\frac{\log p}{\log C_E}}

\DeclareMathOperator*{\Van}{Van}

\newcommand{\Th}{\Theta}

\newcommand{\pab}{P^{(a,b)}}

\newcommand{\Km}{K^{(m)}}
\newcommand{\km}{k^{(m)}}
\newcommand{\Ba}{B^{(a)}}

\newcommand{\half}{{\frac12}}
\newcommand{\mhalf}{{-\frac12}}

\newcommand{\ve}{\varepsilon}

\title{Investigations of Zeros Near the Central Point of Elliptic Curve
$L$-Functions}

\author{Steven J. Miller}
\address{Department of Mathematics, Brown University, 151 Thayer
 Street,
Providence, RI 02912}
 \email{sjmiller@math.brown.edu}

\date{October 10\textsuperscript{th}, 2005}



\begin{document}

\maketitle \centerline{Mathematics Department}
 \centerline{Brown University}
 \centerline{151 Thayer Street}
 \centerline{Providence, RI 02912}

\begin{abstract} We
explore the effect of zeros at the central point on nearby zeros of
elliptic curve $L$-functions, especially for one-parameter families
of rank $r$ over $\Q$. By the Birch and Swinnerton Dyer Conjecture
and Silverman's Specialization Theorem, for $t$ sufficiently large
the $L$-function of each curve $E_t$ in the family has $r$ zeros
(called the family zeros) at the central point. We observe
experimentally a repulsion of the zeros near the central point, and
the repulsion increases with $r$. There is greater repulsion in the
subset of curves of rank $r+2$ than in the subset of curves of rank
$r$ in a rank $r$ family. For curves with comparable conductors, the
behavior of rank $2$ curves in a rank $0$ one-parameter family over
$\Q$ is statistically different from that of rank $2$ curves from a
rank $2$ family. Unlike excess rank calculations, the repulsion
decreases markedly as the conductors increase, and we conjecture
that the $r$ family zeros do not repel in the limit. Finally, the
differences between adjacent normalized zeros near the central point
are statistically independent of the repulsion, family rank and rank
of the curves in the subset. Specifically, the differences between
adjacent normalized zeros are statistically equal for all curves
investigated with rank $0$, $2$ or $4$ and comparable conductors
from one-parameter families of rank $0$ or $2$ over $\Q$.
\end{abstract}


\section{Introduction}

Random matrix theory has successfully modeled the behavior of the
zeros and values of many $L$-functions; see for example the
excellent surveys \cite{KeSn2,Far}. The correspondence first
appeared in Montgomery's analysis of the pair
correlation\footnote{If $\{\alpha_j\}_{j=1}^\infty$ is an increasing
sequence of numbers and $B \subset \R^{n-1}$ is a compact box, the
$n$-level correlations are
\begin{eqnarray}\label{eqnlevelcorr}
\lim_{N \rightarrow \infty}  \# \frac{\{
(\alpha_{j_1}-\alpha_{j_2}, \dots, \alpha_{j_{n-1}} -
\alpha_{j_n}) \in B, j_i \le N, j_i \neq j_k \} }{N} \nonumber\
\end{eqnarray}
One may replace the boxes with smooth test functions; see
\cite{RuSa} for details.} of the zeros of the Riemann zeta
function as the zeros tend to infinity \cite{Mon}. Dyson noticed
that Montgomery's answer, though limited to test functions
satisfying certain support restrictions, agrees with the pair
correlation of the eigenvalues from the Gaussian Unitary
Ensemble\footnote{The GUE is the $N\to\infty$ scaling limit of
$N\times N$ complex Hermitian matrices with entries independently
chosen from Gaussians; see \cite{Meh} for details.} (GUE).
Montgomery conjectured that his result holds for all correlations
and all support. Again with suitable restrictions and in the limit
as the zeros tend to infinity, Hejhal \cite{Hej} showed the triple
correlation of zeros of the Riemann zeta function agree with the
GUE, and, more generally, Rudnick and Sarnak \cite{RuSa} showed
the $n$-level correlations of the zeros of any principal
$L$-function (the $L$-function attached to a cuspidal automorphic
representation of ${\rm GL}_m$ over $\Q$) also agree with the GUE.

In this paper we explore another connection between $L$-functions
and random matrix theory, the effect of multiple zeros at the
central point on nearby zeros of an $L$-function and the effect of
multiple eigenvalues at $1$ on nearby eigenvalues in a classical
compact group. Particularly interesting cases are families of
elliptic curve $L$-functions. It is conjectured that zeros of
primitive $L$-functions are simple, except potentially at the
central point for arithmetic reasons. For an elliptic curve $E$, the
Birch and Swinnerton-Dyer Conjecture \cite{BS-D1,BS-D2} states that
the rank of the Mordell-Weil group $E(\Q)$ equals the order of
vanishing of the $L$-function $L(E,s)$ at the central
point\footnote{We normalize all $L$-functions to have functional
equation $s \mapsto 1-s$, and thus central point is at $s=\foh$.}
$s=\foh$. Let $\mathcal{E}$ be a one-parameter family of elliptic
curves over $\Q$ with (Mordell-Weil) rank\footnote{The group of
rational function solutions $(x(T), y(T)) \in\Q(T)^2$ to $y^2 = x^3
+ A(T)x + B(T)$ is isomorphic to $\Z^r \oplus \T$, where $\T$ is the
torsion part and $r$ is the rank.} $r$: \be y^2 \ = \ x^3 + A(T) x +
B(T), \ \ A(T), B(T) \in \Z[T]. \ee For all $t$ sufficiently large
each curve $E_t$ in the family $\mathcal{E}$ has rank at least $r$,
by Silverman's specialization theorem \cite{Si2}. Thus we expect
each curve's $L$-function to have at least $r$ zeros at the central
point. We call the $r$ conjectured zeros from the Birch and
Swinnerton-Dyer Conjecture the \emph{family zeros}. Thus, at least
conjecturally, these families of elliptic curves offer an exciting
and accessible laboratory where we can explore the effect of
multiple zeros on nearby zeros.

The main tool for studying zeros near the central point (the
\emph{low-lying zeros}) in a family is the $n$-level density. Let
$\phi(x) = \prod_{i=1}^n \phi_i(x_i)$ where the $\phi_i$ are even
Schwartz functions whose Fourier transforms $\hphi_i$ are compactly
supported. Following Iwaniec-Luo-Sarnak \cite{ILS}, we define the
$n$-level density for the zeros of an $L$-function $L(s,f)$ by
\begin{eqnarray}\label{eq:expformexp}
D_{n, f}(\phi) \ &=& \ \sum_{j_1, \dots, j_n \atop j_k \neq \pm
j_\ell} \phi_1\left(\gamma_{f,j_1}\frac{\log C_f}{2\pi}\right)
\cdots \phi_n\left(\gamma_{f,j_n}\frac{\log C_f}{2\pi}\right);
\end{eqnarray}
$C_f$ is the analytic conductor of $L(s,f)$, whose non-trivial zeros
are $\foh + \i\gamma_{f,j}$. Under GRH, the non-trivial zeros all
lie on the critical line $\Re(s) = \foh$, and thus
$\gamma_{f,j}\in\R$. As $\phi_i$ is Schwartz, note that most of the
contribution is from zeros near the central point. The analytic
conductor of an $L$-function normalizes the non-trivial zeros of the
$L$-function so that, near the central point, the average spacing
between normalized zeros is $1$; it is determined by analyzing the
$\Gamma$-factors in the functional equation of the $L$-function (see
for example \cite{ILS}). For elliptic curves the analytic conductor
is the conductor of the elliptic curve (the level of the
corresponding weight $2$ cuspidal newform from the Modularity
Theorem of \cite{Wi,TaWi,BCDT}).

We order a family $\mathcal{F}$ of $L$-functions by analytic
conductors. Let $\mathcal{F}_N = \{f\in \mathcal{F}: C_f \le N\}$.
The $n$-level density for the family $\mathcal{F}$ with test
function $\phi$ is \be D_{n,\mathcal{F}}(\phi) \ = \
\lim_{N\to\infty} D_{n,\mathcal{F}_N}(\phi), \ee where
\be\label{eq:nlddnfn} D_{n,\mathcal{F}_N} \ = \
\frac1{|\mathcal{F}_N|} \sum_{f\in\mathcal{F}_N} D_{n,f}(\phi). \ee
We can of course investigate other subsets. Other common choices are
$\{f: C_f \in [N,2N]\}$, or, for a one-parameter family
$\mathcal{E}$ of elliptic curves over $\Q$, $\{E_t \in \mathcal{E}:
t\in [N,2N]\}$.

Let ${\rm U}(N)$ be the ensemble of $N\times N$ unitary matrices
with Haar measure. The classical compact groups are sub-ensembles
$G(N)$ of ${\rm U}(N)$; the most frequently encountered ones are
${\rm USp}(2M)$, ${\rm SO}(2N)$ and ${\rm SO}(2N+1)$. Katz and
Sarnak's Density Conjecture \cite{KaSa1,KaSa2} states that as the
conductors tend to infinity, the behavior of the normalized zeros
near the central point equals the $N\to\infty$ scaling limit of the
normalized eigenvalues near $1$ of a classical compact group; see
\eqref{eq:nldslqwert} for an exact statement. In the function field
case, the corresponding classical compact group can be identified
from the monodromy group; in the number field case, however, the
reason behind the identification is often a mystery (see \cite{DM}).
As the eigenvalues of a unitary matrix are of the form
$e^{i\theta}$, we often talk about the eigenangles $\theta$ instead
of the eigenvalues $e^{i\theta}$, and the eigenangle $0$ corresponds
to the eigenvalue $1$.

Using the explicit formula we replace the sums over zeros in
\eqref{eq:expformexp} with sums over the Fourier coefficients at
prime powers. For example, if $E: y^2 = x^3 + Ax + B$ is an elliptic
curve, assuming GRH the non-trivial zeros of the associated
$L$-function \be L(E,s) \ = \ \sum_{n=1}^\infty \gl_E(n) n^{-s} \ee
(normalized to have functional equation $s \mapsto 1-s$) are $\foh +
i\gamma$, $\gamma \in \R$. If $\phi$ is a Schwartz test function,
then the explicit formula for $L(E,s)$ is \bea \sum_{\gamma_j}
\phi\left(\gamma_j {\log C_E \over 2\pi}\right) & = &
\widehat{\phi}(0) + \phi(0) - 2 \sum_p \plogce \ \widehat{\phi}
\left( \plogce \right) \frac{\gl_E(p)}{\sqrt{p}} \nonumber\\ & & -\
2\ \sum_p \plogce\ \widehat{\phi}
\left(\frac{2 \log p}{\log C_E} \right) \frac{\gl_E^2(p)}{p} \nonumber\\
& &\ +\ O\left(\frac{\log \log C_E}{\log C_E}\right); \eea see for
example \cite{Mes,Mil1}. By using appropriate averaging formulas and
combinatorics, the resulting prime power sums in the $n$-level
densities can be evaluated for $\hphi_i$ of suitably restricted
support. The Density Conjecture is that to each family of
$L$-functions $\mathcal{F}$, for any Schwartz test function
$\phi:\R^n\to\R^n$, \bea\label{eq:nldslqwert}
D_{n,\mathcal{F}}(\phi) \ = \ \int \phi(x) W_{n,\mathcal{G}}(x)dx \
= \ \int \hphi(u) \widehat{W}_{n,\mathcal{G}}(u)du. \eea The density
kernel $W_{n,\mathcal{G}}(x)$ is determined from the $N\to\infty$
scaling limit of the associated classical compact group $G(N)$; the
last equality follows by Plancherel. The most frequently occurring
answers are the scaling limits of Unitary, Symplectic and Orthogonal
ensembles. For $n=1$ we have
\begin{align}\label{eq:1ldallsupp}
\begin{array}{lcl}
\widehat{W}_{1,{\rm U} }(u) & = & \delta(u) \\
\widehat{W}_{1,\sy}(u) & = & \delta(u) - \foh {\rm I}(u) \\
\widehat{W}_{1,\soe }(u) &\ =\ & \delta(u) + \foh {\rm I}(u) \\
\widehat{W}_{1,\soo }(u) & = & \delta(u) - \foh {\rm I}(u) + 1 \\
\widehat{W}_{1,\so }(u) & = & \delta(u) + \foh,
\end{array}\end{align}
where ${\rm I}(u)$ is the characteristic function of $[-1,1]$. For
arbitrarily small support, unitary and symplectic are
distinguishable from each other and the orthogonal groups; however,
for test functions $\hphi$ supported in $(-1,1)$, the three
orthogonal groups agree: \begin{align}\label{eq:FTphi1ld}
\begin{array}{lcl}
\int \hphi(u)\widehat{W}_{1,\un}(u)du & = & \hphi(u) \\ \int
\hphi(u)\widehat{W}_{1,\sy}(u)du & = & \hphi(u) - \foh \phi(0)\\
\int \hphi(u) \widehat{W}_{1,\soe}(u)du & = & \hphi(u) + \foh \phi(0) \\
\int \hphi(u)\widehat{W}_{1,\soo}(u)du &\ =\ & \hphi(u) + \foh
\phi(0) \\ \int \hphi(u)\widehat{W}_{1,\so}(u)du & = & \hphi(u) +
\foh\phi(0).
 \end{array}\end{align}
Similar results hold for the $n$-level densities, but below we
only need the $1$-level; see \cite{Con,KaSa1} for the derivations
of the general $n$-level densities, and Appendix \ref{app:eduardo}
for the $1$-level density for the orthogonal groups.

For one-parameter families of elliptic curves, the results suggest
that the correct models are orthogonal groups (if all functional
equations are even then the answer is $\soe$, while if all are odd
the answer is $\soo$). Often instead of normalizing each curve's
zeros by the logarithm of its conductor (the local rescaling), one
instead uses the average log-conductor (the global rescaling). If we
are only interested in the average rank, it suffices to study just
the $1$-level density from the global rescaling. This is because we
only care about the imaginary parts of the zeros at the central
point, and both scalings of the imaginary part of the central point
are zero; see for example \cite{Br,Go2,H-B,Mic,Si3,Yo2}. To date all
results have support in $(-1,1)$, where \eqref{eq:FTphi1ld} shows
that the behavior of $\so$, $\soe$ and $\soo$ are indistinguishable.
If we want to specify a unique corresponding classical compact group
we study the $2$-level density as well, which for arbitrarily small
support suffices to distinguish the three orthogonal candidates.
Using the global rescaling removes many complications in the
$1$-level sums but not in the $2$-level sums. In fact, for the
$2$-level investigations the global rescaling is as difficult as the
local rescaling; see \cite{Mil2} for details.

Our research was motivated by investigations on the distribution
of rank in families of elliptic curves as the conductors grow. As
we see below, for the ranges of conductors studied there is poor
agreement between elliptic curve rank data and the $N\to\infty$
scaling limits of random matrix theory. The purpose of this
research is to show that another statistic, the distribution of
the first few zeros above the central point, converges more
rapidly.

We briefly review the excess rank phenomenon. A generic
one-parameter family of elliptic curves over $\Q$ has half of its
functional equations even and half odd (see \cite{He} for the
precise conditions for a family to be generic). Consider such a
one-parameter family of elliptic curves over $\Q$, of rank $r$, and
assume the Birch and Swinnerton-Dyer Conjecture. It is believed that
the behavior of the non-family zeros is modeled by the $N\to\infty$
scaling limit of orthogonal matrices. Thus if the Density Conjecture
is correct, then at the central point in the limit as the conductors
tend to infinity the $L$-functions have exactly $r$ zeros $50\%$ of
the time, and exactly $r+1$ zeros $50\%$ of the time. Thus in the
limit half the curves have rank $r$ and half have rank $r+1$. In a
variety of families, however, one observes\footnote{Actually, this
is not quite true. The analytic rank is estimated by the location of
the first non-zero term in the series expansion of $L(E,s)$ at the
central point (see \cite{Cr} for the algorithms). For example, if
the zeroth through third coefficients are smaller than $10^{-5}$ and
the fourth is $1.701$, then we say the curve has analytic rank $4$,
even though it is possible (though unlikely) that one of the first
four coefficients is really non-zero. It is difficult to prove an
elliptic curve $L$-function vanishes to order two or greater.
Goldfeld \cite{Go1} and Gross-Zagier \cite{GZ} give an effective
lower bound for the class number of imaginary quadratic fields by an
analysis of an elliptic curve $L$-function which is proven to have
three zeros at the central point.}\label{foot:estrank} that $30\%$
to $40\%$ have rank $r$, about $48\%$ have rank $r+1$, $10\%$ to
$20\%$ have rank $r+2$, and about $2\%$ have rank $r+3$; see for
example \cite{BM,Fe1,Fe2,ZK}.

We give a representative family below; see in particular
\cite{Fe2} for more examples. Consider the one-parameter family
$y^2 = x^3 + 16Tx + 32$ of rank $0$ over $\Q$. Each range below
has $2000$ curves:

\begin{center}
\begin{tabular}{rrrrrr} \\
\underline{\ \ \ \ \ \ \ \ \ \ \ $T$-range} & \underline{rank 0} &
\underline{rank 1} &
\underline{rank 2} & \underline{rank 3}  & \underline{run time (hours)}   \\
$[-1000,\ \ 1000)$ & $39.4\%$ & $47.8\%$ & $12.3\%$ & $0.6\%$ & <1 \\

$[1000,\ \ 3000)$ & $38.4\%$ & 47.3\% & 13.6\% & 0.6\% & <1 \\

$[4000,\ \ 6000)$ & $37.4\%$ & 47.8\% & 13.7\% & 1.1\% & 1 \\

$[8000,\ 10000)$ & $37.3\%$ & 48.8\% & 12.9\% & 1.0\% & 2.5 \\

$[24000,26000)$ & $35.1\%$ & 50.1\% & 13.9\% & 0.8\% & 6.8 \\

$[50000,52000)$ & $36.7\%$ & 48.3\% & 13.8\% & 1.2\% & 51.8 \\

\end{tabular} \end{center}

The relative stability of the percentage of curves in a family
with rank $2$ above the family rank $r$ naturally leads to the
question as to whether or not this persists in the limit; it
cannot persist if the Density Conjecture (with orthogonal groups)
is true for all support\footnote{Explicitly, if the
large-conductor limit of the elliptic curve $L$-functions agree
with the $N\to\infty$ scaling limits of orthogonal groups.}.
Recently Watkins \cite{Wat} investigated the family $x^3 + y^3 =
m$ for varying $m$, and unlike other families his range of $m$ was
large enough to see the percentage with rank $r+2$ markedly
decrease, providing support for the Density Conjecture (with
orthogonal groups).

In our example above, as well as the other families investigated,
the logarithms of the conductors are quite small. Even in our last
set the log-conductors are under $40$. An analysis of the error
terms in the explicit formula suggests the rate of convergence of
quantities related to zeros of elliptic curves is like the logarithm
of the conductors. It is quite satisfying when we study the first
few normalized zeros above the central point that, unlike excess
rank, we see a dramatic decrease in repulsion with modest increases
in conductor.

In \S\ref{sec:twoRMTmodels} we study two random matrix ensembles
which are natural candidates to model families of elliptic curves
with positive rank. Many natural questions concerning the
normalized eigenvalues for these models for finite $N$ lead to
quantities that are expressed in terms of eigenvalues of integral
equations. Our hope is that showing the possible connections
between these models and number theory will spur interest in
studying these models and analyzing these integral equations. We
assume the Birch and Swinnerton-Dyer Conjecture, as well as GRH.
We calculate some properties of these ensembles in Appendix
\ref{sec:hard-edge-ensembl}.

In \S\ref{sec:theoryresults} we summarize the theoretical results
of previous investigations, which state:

\bi \item For one-parameter families of rank $r$ over $\Q$ and
suitably restricted test functions, as the conductors tend to
infinity the $1$-level densities imply that in this restricted
range, the $r$ family zeros at the central point are independent
of the remaining zeros. \ei

If this were to hold for all test functions, then as the
conductors tend to infinity the distribution of the first few
normalized zeros above the central point would be independent of
the $r$ family zeros.

In \S\ref{sec:ExpResults} we numerically investigate the first few
normalized zeros above the central point for elliptic curves from
many families of different rank.   Our main observations are:

\bi \item The first few normalized zeros are repelled from the
central point. The repulsion increases with the number of zeros at
the central point, and even in the case when there are no zeros at
the central point there is repulsion from the large-conductor
limit theoretical prediction. This is observed for the family of
all elliptic curves, and for one-parameter families of rank $r$
over $\Q$.

\item There is \emph{greater} repulsion in the first normalized
zero above the central point for subsets of curves of rank $2$
from one-parameter families of rank $0$ over $\Q$ than for subsets
of curves of rank $2$ from one-parameter families of rank $2$ over
$\Q$. It is conjectured that as the conductors tend to infinity,
0\% of curves in a family of rank $r$ have rank $r+2$ or greater.
If this is true, we are comparing a subset of zero relative
measure to one of positive measure. As the first set is
(conjecturally) so small, it is not surprising that to date there
is no known theoretical agreement with any random matrix model for
this case.

\item Unlike most excess rank investigations, as the conductors
increase the repulsion of the first few normalized zeros markedly
decreases. This supports the conjecture that, in the limit as the
conductors tend to infinity, the family zeros are independent of
the remaining normalized zeros (i.e., the repulsion from the
family zeros vanishes in the limit).

\item The repulsion from additional zeros at the central point
cannot entirely be explained by collapsing some zeros to the central
point and leaving all the other zeros alone. See in particular
Remark \ref{rek:attraction}.

\item While the first few normalized zeros are repelled from the
central point, the \emph{differences} between normalized zeros near
the central point are statistically independent of the repulsion, as
well as the method of construction. Specifically, the differences
between adjacent zeros near the central point from curves of rank
$0$, $2$ or $4$ with comparable conductors from one-parameter
families of rank $0$ or $2$ over $\Q$ are statistically equal. Thus
the data suggests that the effect of the repulsion is simply to
shift all zeros by approximately the same amount.

\ei

The numerical data is similar to excess rank investigations. While
both seem to contradict the Density Conjecture, the Density
Conjecture describes the limiting behavior as the conductors tend
to infinity. The rate of convergence is expected to be on the
order of the logarithms of the conductors, which is under $40$ for
our curves. Thus our experimental results are likely misleading as
to the limiting behavior. It is quite interesting that, unlike
most excess rank investigations, we can easily go far enough to
see conductor dependent behavior.

Thus our theoretical and numerical results, as well as the Birch
and Swinnerton-Dyer and Density Conjectures, lead us to

\begin{conj} Consider one-parameter families of elliptic curves of rank
$r$ over $\Q$ and their sub-families of curves with rank exactly
$r+k$ for $k \in \{0,1,2,\dots\}$. For each sub-family there are
$r$ family zeros at the central point, and these zeros repel the
nearby normalized zeros. The repulsion increases with $r$ and
decreases to zero as the conductors tend to infinity, implying
that in the limit the $r$ family zeros are independent of the
remaining zeros. If $k \ge 2$ these additional non-family zeros at
the central point may influence nearby zeros, even in the limit as
the conductors tend to infinity. The spacings between adjacent
normalized zeros above the central point are independent of the
repulsion; in particular, it does not depend on $r$ or $k$, but
only on the conductors.
\end{conj}


\section{Random Matrix Models for Families of Elliptic
Curves}\label{sec:twoRMTmodels}

We want a random matrix model for the behavior of zeros from
families of elliptic curve $L$-functions with a prescribed number of
zeros at the central point. We concentrate on models for either one
or two-parameter families over $\Q$, and refer the reader to
\cite{Far} for more on random matrix models. Both of these families
are expected to have orthogonal symmetries. Many people (see for
example \cite{DFK,Go2,GM,Mai,RuSi,Rub2,ST}) have studied families
constructed by twisting a fixed elliptic curve by characters. The
general belief is that such twisting should lead to unitary or
symplectic families, depending on the orders of the characters.

There are two natural models for the corresponding situation in
random matrix theory of a prescribed number of eigenvalues at $1$
in sub-ensembles of orthogonal groups. For ease of presentation we
consider the case of an even number of eigenvalues at $1$; the odd
case is handled similarly.

Consider a matrix in ${\rm SO}(2N)$. It has $2N$ eigenvalues in
pairs $e^{\pm \i \theta_j}$, with $\theta_j \in [0,\pi]$. The joint
probability measure on $\Theta = (\theta_1,\dots,\theta_N) \in
[0,\pi]^N$ is
\begin{align}
  \label{eq:7}
  d\epsilon_0(\Theta) &\ = \ c_N \prod_{1\le j<k \le N}(\cos\theta_k-\cos\theta_j)^2
  \prod_{1\le j\le N} d\theta_j,
  \end{align}
where $c_N$ is a normalization constant so that
$d\epsilon_0(\Theta)$ integrates to $1$. From~\eqref{eq:7} we can
derive all quantities of interest on the random matrix side; in
particular, $n$-level densities, distribution of first normalized
eigenvalue above $1$ (or eigenangle above $0$), and so forth.

We now consider two models for sub-ensembles of ${\rm SO}(2N)$
with $2r$ eigenvalues at 1, and the $N\to\infty$ scaling limit of each.\\

\textbf{Independent Model:} The sub-ensemble of ${\rm SO}(2N)$
with the upper left block a $2r\times 2r$ identity matrix. The
joint probability density of the remaining $N-r$ pairs is given by
\begin{equation}
  \label{eq:14a}
  d\varepsilon_{2r,{\rm Indep}}(\Theta) \ = \ c_{2r,{\rm Indep},N}
  \prod_{1 \le j<k \le N-r}
  (\cos\theta_k-\cos\theta_j)^2
  \prod_{1\le j \le N-r} d\theta_j.
\end{equation} Thus the ensemble is matrices of the form
\begin{equation}
  \label{eq:28}
  \left\{
    \begin{pmatrix}
      I_{2r\times 2r} & \\
      & g
    \end{pmatrix}
    : g \in {\rm SO}(2N-2r)\right\};
\end{equation} the probabilities are equivalent to choosing $g$
with respect to Haar measure on ${\rm SO}(2N-2r)$. We call this the
Independent Model as the forced eigenvalues at 1 from the
$I_{2r\times 2r}$ block do not interact with the eigenvalues of $g$.
In particular, the distribution of the remaining $N-r$ pairs of
eigenvalues is exactly that of ${\rm SO}(2N-2r)$; this block's
$N\to\infty$ scaling limit is just $\soe$. See \cite{Con,KaSa1} as
well
as Appendix \ref{app:eduardo}. \\

\textbf{Interaction Model:} The sub-ensemble of ${\rm SO}(2N)$ with
$2r$ of the $2N$ eigenvalues equaling $1$:
\begin{equation}
  \label{eq:14}
  d\varepsilon_{2r,{\rm Inter}}(\Theta) \ = \ c_{2r,{\rm Inter},N}
  \prod_{1\le j<k\le N-r}
  (\cos\theta_k-\cos\theta_j)^2\prod_{1\le j\le N-r}(1-\cos\theta_j)^{2r}
  d\theta_j.
\end{equation}
We call this the Interaction Model as the forced eigenvalues at 1
\emph{do} affect the behavior of the other eigenvalues near 1.
Note here we condition on all ${\rm SO}(2N)$ matrices with at
least $2r$ eigenvalues equal to $1$. The $(1-\cos\theta_j)^{2r}$
factor results in the forced eigenvalues at $1$ repelling the
nearby eigenvalues.


\begin{rek} As the calculations for the local statistics near the
eigenvalue at~$1$ in the Interaction Model has not appeared in
print, in Appendix~\ref{app:eduardo} (written by Eduardo Due\~nez)
is a derivation of formula~\eqref{eq:14} (see especially
\S\ref{sec:Hard-SO}), as well as the relevant integral (Bessel)
kernels dictating such statistics. See also \cite{Sn} for the
value distribution of the first non-zero derivative of the
characteristic polynomials of this ensemble.
\end{rek}

While both models have at least $2r$ eigenvalues equal to $1$, they
are very different sub-ensembles of ${\rm SO}(2N)$, and they have
distinct limiting behavior (see also Remark
\ref{rek:subensmblesbeh}). We can see this by computing the
$1$-level density for each, and comparing with \eqref{eq:FTphi1ld}.
Letting $\widehat{W}_{1,\soe}$ (respectively
$\widehat{W}_{1,\soe,{\rm Indep},2r}$ and $\widehat{W}_{1,\soe,{\rm
Inter},2r}$) denote the Fourier transform of the kernel for the
$1$-level density of $\soe$ (respectively, of the Independent Model
for the sub-ensemble of $\soe$ with $2r$ eigenvalues at $1$ and of
the Interaction Model for the sub-ensemble of $\soe$ with $2r$
eigenvalues at $1$), we find in Appendix \ref{app:eduardo} that
\bea\label{eq:three1LDs}
  \widehat{W}_{1,\soe}(u)& \ =\ &  \delta(u) + \foh {\rm I}(u) \nonumber\\
  \widehat{W}_{1,\soe,{\rm Indep},2r}(u) &  = &  \delta(u) + \foh {\rm I}(u) +
  2\nonumber\\
  \widehat{W}_{1,\soe,{\rm Inter},2r}(u)& =&
  \delta(u) + \foh {\rm I}(u) +  2 + 2(|u|-1)
  {\rm I}(u).\eea
As ${\rm I}(u)$ is positive for $|u| < 1$, note that the density is
smaller for $|u| < 1$ in the Interaction versus the Independent
Model. We can interpret this as a repulsion of zeros, as the
following heuristic shows (though see Appendix \ref{app:eduardo} for
proofs). We compare the 1-level density of zeros from curves with
and without repulsion, and show that for a positive decreasing test
function, the 1-level density is smaller when there is repulsion.

Consider two elliptic curves, $E$ of rank $0$ and conductor $C_E$
and $E'$ of rank $r$ and conductor $C_{E'}$. Assume $C_E \approx
C_{E'} \approx C$, and assume GRH for both $L$-functions. If the
curve $E$ has rank $0$ then we expect the first zero above the
central point, $\foh + \i\gamma_{E,1}$, to have $\gamma_{E,1}
\approx \frac{1}{\log C}$. For $E_r$, if the $r$ family zeros at the
central point repel, it is reasonable to posit a repulsion of size
$\frac{b_r}{\log C}$ for some $b_r>0$. This is because the natural
scale for the distance between the low-lying zeros is $\frac1{\log
C}$, so we are merely positing that the repulsion is proportional to
the distance. We assume all zeros are repelled equally; evidence for
this is provided in \S\ref{sec:spacingsnorzeros}. Thus for $E'$ (the
repulsion case) we assume $\gamma_{E',j} \approx \gamma_{E,j} +
\frac{b_r}{\log C}$. We can detect this repulsion by comparing the
1-level densities of $E$ and $E'$. Take a non-negative decreasing
Schwartz test function $\phi$. The difference between the
contribution from the $j$\textsuperscript{th} zero of each is \bea
\phi\left(\gamma_{E',j}\frac{\log C}{2\pi}\right) -
\phi\left(\gamma_{E,j}\frac{\log C}{2\pi}\right) & \ \approx \ &
\phi\left(\gamma_{E,j}\frac{\log C}{2\pi} + \frac{b_r}{2\pi}\right)
- \phi\left(\gamma_{E,j}\frac{\log C}{2\pi}\right) \nonumber\\ &
\approx & \phi'\left(\gamma_{E,j}\frac{\log C}{2\pi}\right) \cdot
\frac{b_r}{2\pi}. \eea As $\hphi$ is decreasing, its derivative is
negative and thus the above shows the 1-level density for the zeros
from $E'$ (assuming repulsion) is smaller than the 1-level density
for zeros from $E$. Thus the lower $1$-level density in the
Interaction Model versus the Independent Model can be interpreted as
a repulsion; however, this repulsion can be shared among several
zeros near the central point. In fact, the observations in
\S\ref{sec:spacingsnorzeros} suggest that the repulsion shifts all
normalized zeros near the central point approximately equally.


\section{Theoretical Results}\label{sec:theoryresults}

Consider a one-parameter family of elliptic curves of rank $r$
over $\Q$. We summarize previous investigations of the effect of
the (conjectured) $r$ family zeros on the other zeros near the
central point. For convenience we state the results for the global
rescaling, though similar results hold for the local rescaling
(under slightly more restrictive conditions; see \cite{Mil2} for
details). For small support, the $1$ and $2$-level densities agree
with the scaling limits of \be \mattwo{I_{r\times r}}{}{}{{\rm
O}(N)}, \ \ \ \mattwo{I_{r\times r}}{}{}{{\rm SO}(2N)}, \ \ \
\mattwo{I_{r\times r}}{}{}{{\rm SO}(2N+1)}, \ee depending on
whether or not the signs of the functional equation are
equidistributed or all the signs are even or all the signs are
odd. The $1$ and $2$-level densities provide evidence towards the
Katz-Sarnak Density Conjecture for test functions whose Fourier
transforms have small support (the support is computable and
depends on the family). See \cite{Mil1} for the calculations with
the global rescaling, though the result for the $1$-level density
is implicit in \cite{Si3}. Similar results are observed for
two-parameter families of elliptic curves in \cite{Mil1,Yo2}.

While the above results are consistent with the Birch and
Swinnerton-Dyer Conjecture that each curve's $L$-function has at
least $r$ zeros at the central point, it is not a proof (even in the
limit) because our supports are finite. For families with $t \in
[N,2N]$ the errors are of size $O(\frac{1}{\log N})$ or
$O(\frac{\log\log N}{\log N})$. Thus for large $N$ we cannot
distinguish a family with exactly $r$ zeros at the central point
from a family where each $E_t$ has exactly $r$ zeros at $\pm(\log
C_t)^{-2007}$.

For one-parameter families of elliptic curves over $\Q$, in the
limit as the conductors tend to infinity the family zeros (those
arising from our belief in the Birch and Swinnerton-Dyer Conjecture)
appear to be independent from the other zeros. Equivalently, if we
remove the contributions from the $r$ family zeros, for test
functions with suitably restricted support the spacing statistics of
the remaining zeros agree perfectly with the standard orthogonal
groups $\so$, $\soe$ and $\soo$, and it is conjectured that these
results should hold for all support. Thus the $n$-level density
arguments support the Independent over the Interaction Model when we
study \emph{all} curves in a family; however, these theoretical
arguments do not apply if we study the sub-family of curves of rank
$r+k$ ($k \ge 2$) in a rank $r$ one-parameter family over $\Q$.

\begin{rek}\label{rek:subensmblesbeh}
It is important to note that our theoretical results are for the
entire one-parameter family. Specifically, consider the subset of
curves of rank $r+2$ from a one-parameter family of rank $r$ over
$\Q$. If the Density Conjecture (with orthogonal groups) is true,
then in the limit 0\% of curves are in this sub-family. Thus these
curves may behave differently without contradicting the
theoretical results for the entire family. Situations where
sub-ensembles behave differently than the entire ensemble are well
known in random matrix theory. For example, to any simple graph we
may attach a real symmetric matrix, its adjacency matrix, where
$a_{ij} = 1$ if there is an edge connecting vertices $i$ and $j$,
and $0$ otherwise. The adjacency matrices of $d$-regular graphs
are a thin sub-ensemble of real symmetric matrices with entries
independently chosen from $\{-1,0,1\}$. The density of normalized
eigenvalues in the two cases are quite different, given by
Kesten's Measure \cite{McK} for $d$-regular graphs and Wigner's
Semi-Circle Law \cite{Meh} for the real symmetric matrices.
\end{rek}

It is an interesting question to determine the appropriate random
matrix model for rank $r+2$ curves in a rank $r$ one-parameter
family over $\Q$, both in the limit of large conductors as well as
for finite conductors. We explore this issue in greater detail in
\S\ref{sec:oneparamrank0} to \S\ref{sec:spacingsnorzeros}, where
we compare the behavior of rank $2$ curves from rank $0$
one-parameter families over $\Q$ to that of rank $2$ curves from
rank $2$ one-parameter families over $\Q$.


\section{Experimental Results}\label{sec:ExpResults}

We investigate the first few normalized zeros above the central
point. We used Michael Rubinstein's $L$-function calculator
\cite{Rub3} to determine the zeros. The program does a contour
integral to ensure that we found all the zeros in a region, which
is essential in studies of the first zero! See \cite{Rub1} for a
description of the algorithms. The analytic ranks were found (see
Footnote \ref{foot:estrank}) by determining the values of the
$L$-functions and their derivatives at the central point by the
standard series expansion; see \cite{Cr} for the algorithms. Some
of the programs and all of the data (minimal model, conductor,
discriminant, sign of the functional equation, first non-zero
Taylor coefficient from the series expansion at the central point,
and the first three zeros above the central point) are available
online at
\begin{center}
\texttt{http://www.math.brown.edu/$\sim$sjmiller/repulsion}
\end{center}

We study several one-parameter families of elliptic curves over
$\Q$. As all of our families are rational surfaces\footnote{An
elliptic surface $y^2 = x^3 + A(T)x + B(T)$, $A(T)$, $B(T) \in
\Z[T]$, is a rational surface if and only if one of the following
is true: $(1)\ $ $0 < \max\{3 \mbox{deg} A, 2\mbox{deg} B\} < 12$;
$(2)\ $ $3\mbox{deg} A = 2\mbox{deg} B = 12$ and
$\mbox{ord}_{T=0}T^{12} \Delta(T^{-1}) = 0$.}, Rosen and
Silverman's result that the weighted average of fibral Frobenius
trace values determines the rank over $\Q$ (see \cite{RoSi}) is
applicable, and evaluating simple Legendre sums suffices to
determine the rank. We mostly use one-parameter families from
Fermigier's tables \cite{Fe2}, though see \cite{ALM} for how to
use the results of \cite{RoSi} to construct additional
one-parameter families with rank over $\Q$.

We cannot obtain a decent number of curves with approximately
equal log-conductors by considering a solitary one-parameter
family. The conductors in a family typically grow polynomially in
$t$. The number of Fourier coefficients needed to study a value of
$L(s,E_t)$ on the critical line is of order $\sqrt{C_t}\log C_t$
($C_t$ is the conductor of $E_t$), and we must then additionally
evaluate numerous special functions. We can readily calculate the
needed quantities up to conductors of size $10^{11}$, which
usually translates to just a few curves in a family. We first
studied all elliptic curves (parametrized with more than one
parameter), found the minimal models, and then sorted by
conductor. We then studied several one-parameter families,
amalgamating data from different families if the curves had the
same rank and similar log-conductor.

\begin{rek}\label{rek:nucphys}
Amalgamating data from different one-parameter families warrants
some discussion. We expect that the behavior of zeros from curves
with similar conductors and the same number of zeros and family
zeros at the central point should be approximately equal. In other
words, we hope that curves with the same rank and approximately
equal conductors from different one-parameter families of the same
rank $r$ over $\Q$ behave similarly, and we may regard the different
one-parameter families of rank $r$ over $\Q$ as different
measurements of this universal behavior. This is similar to
numerical investigations of the spacings of energy levels of heavy
nuclei; see for example \cite{HH,HPB}. In studying the spacings of
these energy levels, there were very few (typically between 10 and
100) levels for each nucleus. The belief is that nuclei with the
same angular momentum and parity should behave similarly. The
resulting amalgamations often have thousands of spacings and
excellent agreement with random matrix predictions.
\end{rek}

Similar to the excess rank phenomenon, we found disagreement
between the experiments and the predicted large-conductor limit;
however, we believe that this disagreement is due to the fact that
the logarithms of the conductors investigated are small. In
\S\ref{sec:exp0} to \S\ref{sec:compareoneparams} we find that for
curves with zeros at the central point, the first normalized zero
above the central point \emph{is} repelled, and the more zeros at
the central point, the greater the repulsion. However, the
repulsion decreases as the conductors increase. Thus the repulsion
is probably due to the small conductors, and in the limit the
Independent Model (which agrees with the function field analogue
and the theoretical results of \S\ref{sec:theoryresults}) should
correctly describe the first normalized zero above the central
point in curves of rank $r$ in families of rank $r$ over $\Q$. It
is not known what the correct model is for curves of rank $r+2$ in
a family of rank $r$ over $\Q$, though it is reasonable to
conjecture it is the Interaction model with the sizes of the
matrices related to the logarithms of the conductors. Keating and
Snaith \cite{KeSn1,KeSn2} showed that to study zeros at height $T$
it is better to look at $N\times N$ matrices, with $N = \log T$,
than to look at the $N\to\infty$ scaling limit. A fascinating
question is to determine the correct finite conductor analogue for
the two different cases here. Interestingly, we see in
\S\ref{sec:spacingsnorzeros} that the spacings between adjacent
normalized zeros is statistically independent of the repulsion,
which implies that the effect of the zeros at the central point
(for finite conductors) is to shift \emph{all} the nearby zeros
approximately equally.


\subsection{Theoretical Predictions: Independent Model}

In Figures $1$ and $2$ we plot the first normalized eigenangle above
$0$ for ${\rm SO}(2N)$ (i.e., $\soe$) and ${\rm SO}(2N+1)$ (i.e.,
$\soo$) matrices. The eigenvalues occur in pairs $e^{\pm \i
\theta_j}$, $\theta_j \in [0,\pi]$; by normalized eigenangles for
$\soe$ or $\soo$ we mean $\theta_j \frac{N}{\pi}$. We chose $2N \le
6$ and $2N+1=7$ for our simulations, and chose our matrices with
respect to the appropriate Haar measure\footnote{Note that for
$\soo$ matrices there is always an eigenvalue at $1$. The
$N\to\infty$ scaling limit of the distribution of the second
eigenangle for $\soo$ matrices equals the $N\to\infty$ scaling limit
of the distribution of the first eigenangle for $\sy$ (Unitary
Symplectic) matrices; see pages 10--11 and 411--416 of \cite{KaSa1}
and page 10 of \cite{KaSa2}.}. We thank Michael Rubinstein for
sharing his $N\to\infty$ scaling limit plots for ${\rm SO}(2N)$ and
${\rm SO}(2N+1)$.


\newpage

\begin{figure}[!h]
\begin{center}
  \includegraphics[width=5cm]{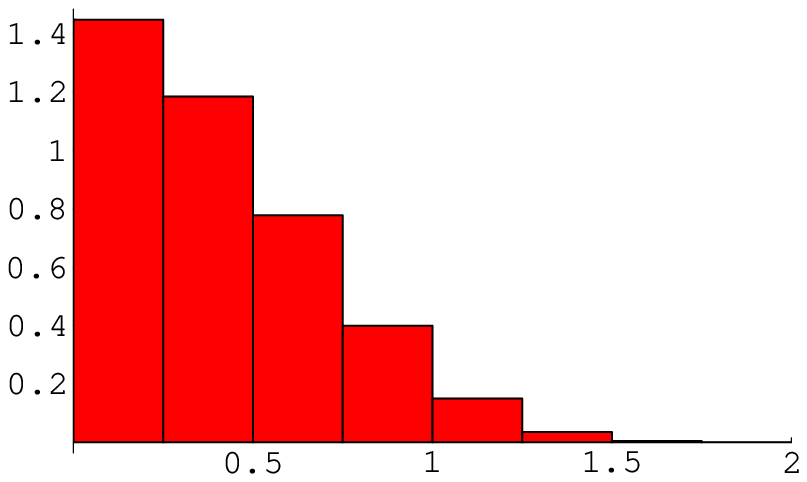}\\
 Figure 1a: First normalized eigenangle above $0$: 23,040 SO(4) matrices\\
 Mean $= .357$, Standard Deviation about the Mean $= .302$, Median $= .357$\\
 \includegraphics[width=5cm]{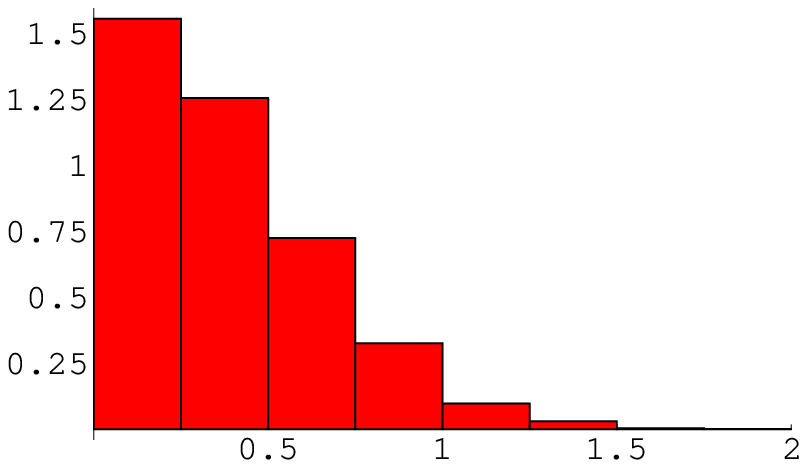}\\
 Figure 1b: First normalized eigenangle above $0$: 23,040 SO(6) matrices\\
 Mean $= .325$, Standard Deviation about the Mean $= .284$, Median $= .325$ \\
  \includegraphics[width=5cm]{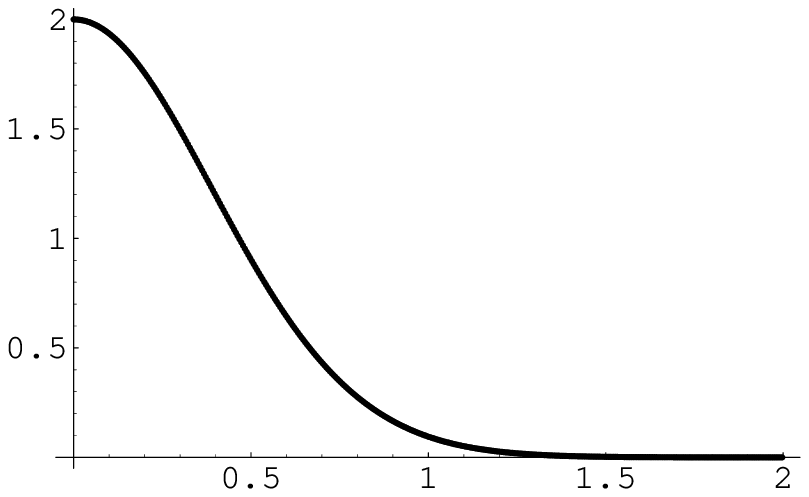}\\
 Figure 1c: First normalized eigenangle above $0$:\\ $N\to\infty$
 scaling limit of SO($2N$): Mean = .321.
 \end{center}
 \end{figure}

\newpage

\begin{figure}
\begin{center}
 \includegraphics[width=5cm]{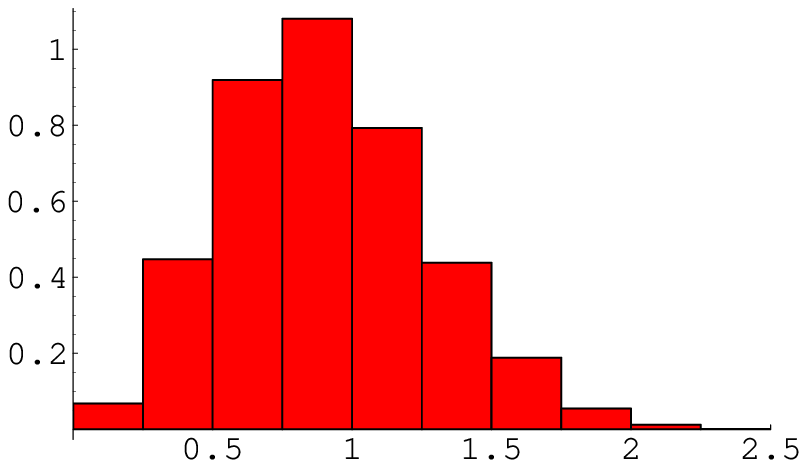}\\
 Figure 2a: First normalized eigenangle above $1$: 322,560 SO(7) matrices\\
 Mean $= .879$, Standard Deviation about the Mean $= .361$, Median $= .879$
 \includegraphics[width=5cm]{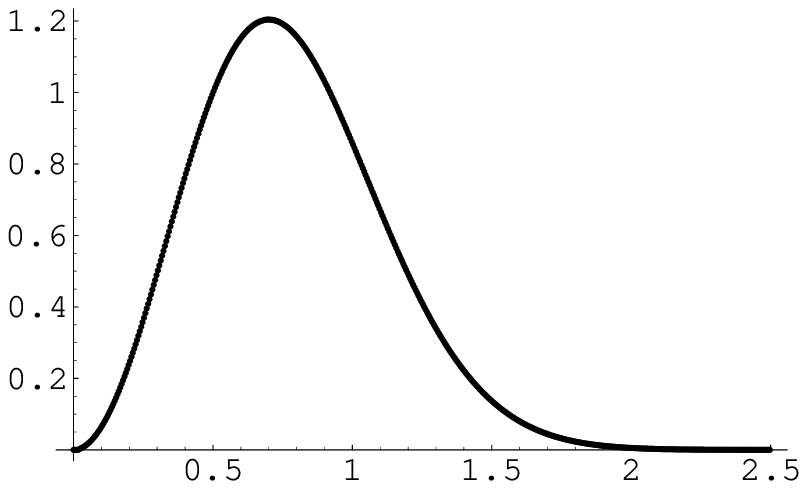}\\
 Figure 2b: First normalized eigenangle above $1$:\\ $N\to\infty$
 scaling limit of SO($2N+1$): Mean = .782.
 \end{center}
 \end{figure}

For the ${\rm SO}(2N)$ matrices, note the mean decreases as $2N$
increases. A similar result holds for ${\rm SO}(2N+1)$ matrices;
as we primarily study even rank below, we concentrate on ${\rm
SO}(2N)$ here. As $N\to\infty$, Katz and Sarnak (pages 412--415 of
\cite{KaSa2}) prove that the mean of the first normalized
eigenangle above $\theta = 0$ (corresponding to the eigenvalue
$1$) for $\soe$ is approximately $0.321$, while for $\soo$ it is
approximately $0.782$.

We study the first normalized zero above the central point for
elliptic curve $L$-functions in \S\ref{sec:exp0} to
\S\ref{sec:compareoneparams}. We rescale each zero: $\gamma_{E_t,1}
\mapsto \gamma_{E_t,1}\frac{\log C_t}{2\pi}$. The mean of the first
normalized eigenangle above $0$ for ${\rm SO}(2N)$ matrices
decreases as $2N$ increases, and similarly we see that the first
normalized zero above the central point in families of elliptic
curves decreases as the conductor increases. This suggests that a
good finite conductor model for families of elliptic curves with
even functional equation and conductors of size $C$ would be ${\rm
SO}(2N)$, with $N$ some function of $\log C$.


\newpage
\subsection{All Curves}\label{sec:exp0}

\subsubsection{Rank $0$ Curves}\label{sec:exp00}

We study the first normalized zero above the central point for
1500 rank $0$ elliptic curves, 750 with $\log({\rm cond}) \in
[3.2,12.6]$ in Figure $3$ and 750 with $\log({\rm cond}) \in
[12.6,14.9]$ in Figure $4$. These curves were obtained as follows:
an elliptic curve can be written in Weierstrass form as \be y^2 +
a_1xy + a_3y \ = \ x^3 + a_2x^2 + a_4x + a_6, \ \ \ a_i \in \Z.
\ee We often denote the curve by $[a_1,a_2,a_3,a_4,a_6]$. We let
$a_1$ range from $0$ to $10$ (as without loss of generality we may
assume $a_1 \ge 0$) and the other $a_i$ range from $-10$ to $10$.
We kept only non-singular curves. We took minimal Weierstrass
models for the ones left, and pruned the list to ensure that all
the remaining curves were distinct. We then analyzed the first few
zeros above the central point for 1500 of these curves (due to the
length of time it takes to compute zeros for the curves, it was
impossible to analyze the entire set).

Figures $3$ and $4$ suggest that as the conductor increases the
repulsion decreases. For the larger conductors in Figure $4$, the
results are closer to the predictions of Katz-Sarnak, and the shape
of the distribution with larger conductors is closer to the random
matrix theory plots of Figure $1$. Though both plots in Figure $3$
and $4$ differ from the random matrix theory plots, the plot in
Figure $4$ is more peaked, the peak occurs earlier, and the decay in
the tail is faster. Standard statistical tests show the two means
(1.04 for the smaller conductors and 0.88 for the larger) are
significantly different. Two possible tests are the Pooled
Two-Sample $t$-Procedure\footnote{\label{foot:pooledt}The Pooled
Two-Sample $t$-Procedure is \be t \ = \
(\overline{X_1}-\overline{X_2}) \Big/ s_p \sqrt{\frac1{n_1}
+\frac1{n_2}}, \ee where $\overline{X_i}$ is the sample mean of
$n_i$ observations of population $i$, $s_i$ is the sample standard
deviation and \be s_p \ = \ \sqrt{\frac{(n_1-1)s_1^2 +
(n_2-1)s_2^2}{n_1+n_2-2}} \ee is the pooled variance; $t$ has a
$t$-distribution with $n_1+n_2-2$ degrees of freedom.} (where we
assume the data are independently drawn from two normal
distributions with the same mean and variance) and the Unpooled
Two-Sample $t$-Procedure\footnote{Notation as in Footnote
\ref{foot:pooledt}, the Unpooled Two-Sample $t$-Procedure is \be t \
= \ (\overline{X_1}-\overline{X_2}) \Big/ \sqrt{\frac{s_1^2}{n_1}
+\frac{s_2^2}{n_2}}; \ee this is approximately a $t$ distribution
with \be \frac{(n_1-1)\ (n_2-1)\ (n_2s_1^2+n_1
s_2^2)^2}{(n_2-1)n_2^2 s_1^4 + (n_1-1)n_1^2 s_2^4}\ee degrees of
freedom} (where we assume the data are independently drawn from two
normal distributions with the same mean and no assumption is made on
the variance). See for example \cite{CaBe}, pages 409-410. Both
tests give $t$-statistics around $10.5$ with over $1400$ degrees of
freedom. As the number of degrees of freedom is so large, we may use
the Central Limit Theorem and replace the $t$-statistic with a
$z$-statistic. As for the standard normal the probability of being
at least $10.5$ standard deviations from zero is less than $3.2
\times 10^{-12}$ percent, we obtain strong evidence against the null
hypothesis that the two means are equal (i.e., we obtain evidence
that the repulsion decreases as the conductor increases).

\begin{figure}[h]
\begin{center}
 \includegraphics[width=7.01cm]{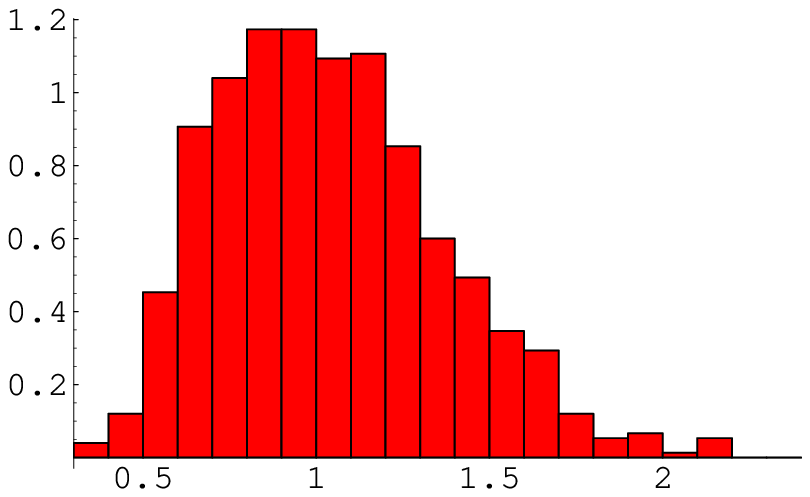}\\
Figure 3: First normalized zero above the central point:\\ $750$
rank 0 curves from $y^2 + a_1xy+a_3y=x^3+a_2x^2+a_4x+a_6$,\\
$\log({\rm cond}) \in [3.2, 12.6]$, $\text{median} = 1.00$,
$\text{mean} = 1.04$,\\ standard deviation about the mean $=.32$
\end{center}
\end{figure}

\begin{figure}[h]
\begin{center}
\includegraphics[width=7.01cm]{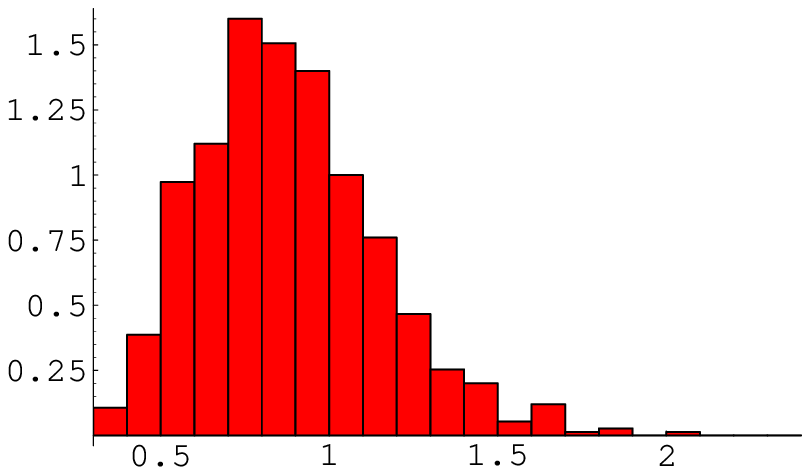}\\
Figure 4: First normalized zero above the central point:\\ $750$
rank 0 curves from $y^2 + a_1xy+a_3y=x^3+a_2x^2+a_4x+a_6$,\\
$\log({\rm cond}) \in [12.6, 14.9]$, $\text{median}=.85$,
$\text{mean} = .88$,\\ standard deviation about the mean $= .27$
\end{center}
\end{figure}

\pagebreak

\newpage
\subsubsection{Rank 2 Curves}\label{sec:exp2}

We study the first normalized zero above the central point for
1330 rank $2$ elliptic curves, 665 with $\log({\rm cond}) \in
[10,10.3125]$ in Figure 5 and 665 with $\log({\rm cond}) \in
[16,16.5]$ in Figure 6. These curves were obtained from the same
procedure which generated the $1500$ curves in \S\ref{sec:exp00},
except now we chose $1330$ curves with what we believe is analytic
rank exactly $2$. We did this by showing the $L$-function has even
sign, the value at the central point is zero to at least 5 digits,
and the second derivative at the central point is non-zero; see
also Footnote \ref{foot:estrank}. In \S\ref{sec:oneparamrank0} and
\S\ref{sec:famrank2oneparamoverQT2} we study other families of
curves of rank $2$ (rank $2$ curves from rank $0$ and rank $2$
one-parameter families over $\Q$).

The results are very noticeable. The first normalized zero is
significantly higher here than for the rank $0$ curves. This
supports the belief that, for small conductors, the repulsion of
the first normalized zero increases with the number of zeros at
the central point.

We again split the data into two sets (Figures 5 and 6) based on the
size of the conductor. As the conductors increase the mean (and
hence the repulsion) significantly decreases, from 2.30 to 1.82.

We are investigating rank $2$ curves from the family of all elliptic
curves (which is a many parameter rank $0$ family). In the limit we
believe half of the curves are rank $0$ and half are rank $1$. The
natural question is to determine the appropriate model for this
subset of curves. As in the limit we believe a curve has rank $2$
(or more) with probability zero, this is a question about
conditional probabilities.

\begin{figure}[!hb]
\begin{center}
\includegraphics[width=7cm]{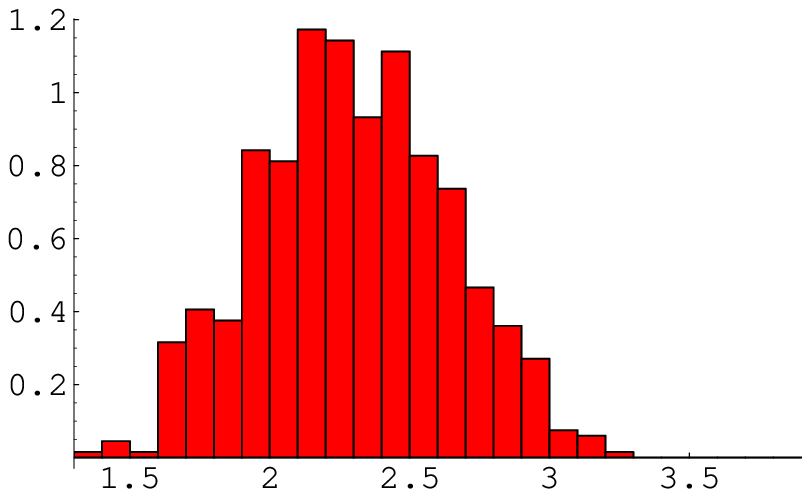}\\
Figure 5: First normalized zero above the central point:\\ $665$
rank $2$ curves from $y^2 + a_1xy+a_3y=x^3+a_2x^2+a_4x+a_6$.\\
$\log({\rm cond}) \in [10, 10.3125]$, $\text{median} = 2.29$,
$\text{mean} = 2.30$
\end{center}
\end{figure}

\begin{figure}[!h]
\begin{center}
\includegraphics[width=7cm]{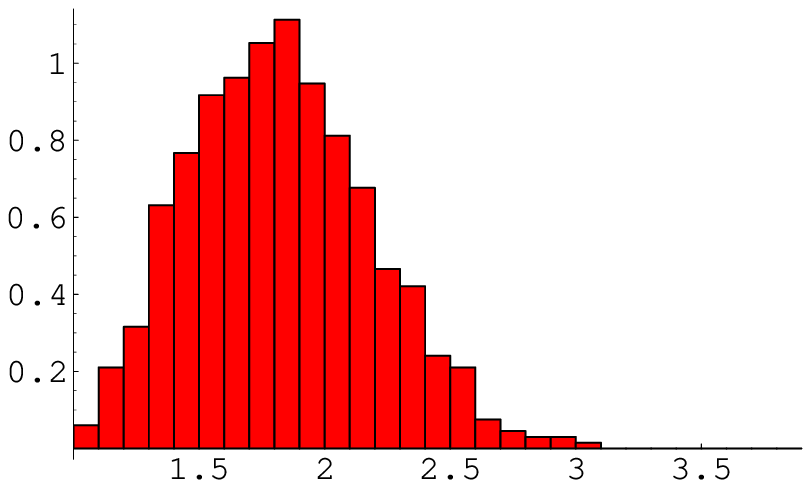}\\
Figure 6: First normalized zero above the central point:\\ $665$
rank $2$ curves from $y^2 + a_1xy+a_3y=x^3+a_2x^2+a_4x+a_6$.\\
$\log({\rm cond}) \in [16, 16.5]$, $\text{median} = 1.81$,
$\text{mean} = 1.82$
\end{center}
\end{figure}


\newpage
\subsection{One-Parameter Families of Rank $0$ Over
$\Q$}\label{sec:oneparamrank0}

\subsubsection{Rank $0$ Curves}

We analyzed $14$ one-parameter families of rank $0$ over $\Q$; we
chose these families from \cite{Fe2}. We want to study rank $0$
curves in a solitary one-parameter family; however, the conductors
grow rapidly and we can only analyze the first few zeros from a
small number of curves in a family. For our conductor ranges it
takes several hours of computer time to find the first few zeros for
all the curves in a family. In Figures 7 and 8 and Tables $1$ and
$2$ we study the first normalized zero above the central point for
$14$ one-parameter families of elliptic curves of rank $0$ over
$\Q$. Even though we have few data points in each family, we note
the medians and means are always higher for the smaller conductors
than the larger ones. Thus the ``repulsion'' is decreasing with
increasing conductor, though perhaps repulsion is the wrong word
here as there is no zero at the central point! We studied the median
as well as the mean because, for small data sets, one or two
outliers can significantly affect the mean; the median is more
robust.

For both the Pooled and Unpooled Two-Sample $t$-Procedure the
$t$-statistic exceeds $20$ with over $200$ degrees of freedom. The
Central Limit Theorem is an excellent approximation and yields a
$z$-statistic exceeding $20$, which strongly argues for rejecting
the null hypothesis that the two means are equal (i.e., providing
evidence that the repulsion decreases with increasing conductors).
Note the first normalized zero above the central point is
significantly larger than the $N\to\infty$ scaling limit of ${\rm
SO}(2N)$ matrices, which is about $0.321$.

Some justification is required for regarding the data from the
$14$ families as independent samples from the same distribution.
It is possible that there are family-specific lower order terms to
the $n$-level densities (see \cite{Mil1,Mil3,Yo2}). Our
amalgamation of the data is similar to physicists combining the
energy level data from different heavy nuclei with similar quantum
numbers. The hope is that the systems are similar enough to
justify such averaging as it is impractical to obtain sufficient
data for just one nucleus (or one family of elliptic curves, as we
see in \S\ref{sec:famrank2oneparamoverQT2}). See also Remark
\ref{rek:nucphys}.

\begin{rek}\label{rek:famnoindep}
The families are not independent: there are $11$ curves that occur
twice and one that occurs three times in the small conductor set
of $220$ curves, and $133$ repeats in the large conductor set of
$996$ curves. In our amalgamations of the families, we present the
results when we double count these curves as well as when we keep
only one curve in each repeated set. In both cases the repeats
account for a sizeable percentage of the total number of
observations; however, there is no significant difference between
the two sets. Any curve can be placed in infinitely many
one-parameter families; given polynomials of sufficiently high
degree we can force any number of curves to lie in two distinct
families. Thus it is not surprising that we run into such problems
when we amalgamate. When we remove the repeated curves, the Pooled
and Unpooled Two-Sample $t$-Procedures still give $t$-statistics
exceeding $20$ with over $200$ degrees of freedom, indicating the
two means significantly differ and supporting the claim that the
repulsion decreases with increasing conductor.
\end{rek}

\begin{figure}
\begin{center}
\includegraphics[width=7cm]{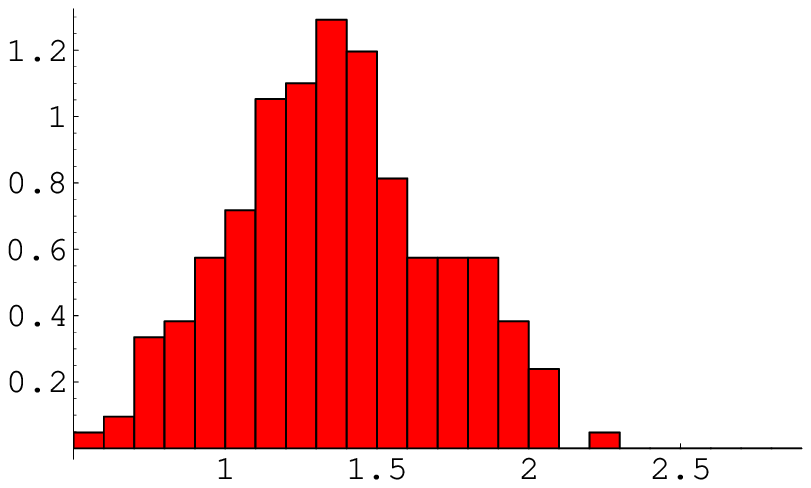}\\ Figure 7: First
normalized zero above the central point.\\ 209 rank 0 curves from
14 rank $0$ one-parameter families,\\ $\log({\rm cond}) \in
[3.26,9.98]$, $\text{median} = 1.35$, $\text{mean} =
1.36$\\ \ \\

\includegraphics[width=7cm]{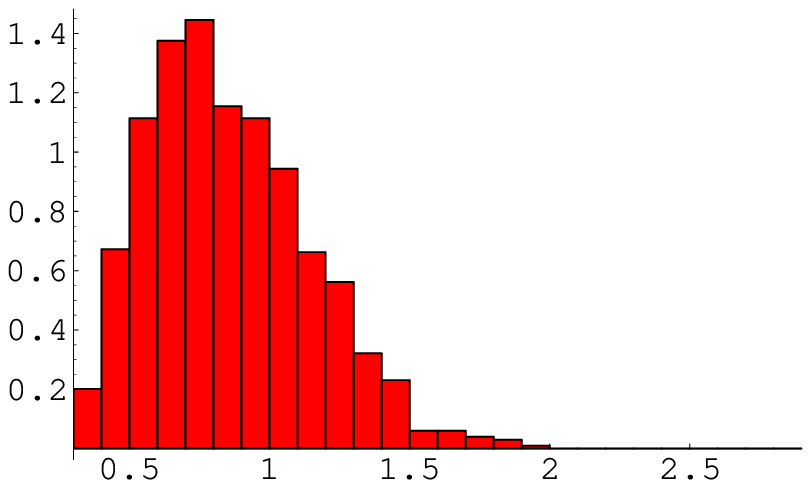}\\ Figure
8: First
normalized zero above the central point.\\
996 rank 0 curves from 14 rank $0$ one-parameter families,\\
$\log({\rm cond}) \in [15.00, 16.00]$, $\text{median} = .81$,
$\text{mean} = .86$.
\end{center}
\end{figure}

\begin{table}[h]
\begin{center}
\caption{First normalized zero above the central point for $14$
one-parameter families of elliptic curves of rank $0$ over $\Q$
(smaller conductors)}
\begin{tabular}{|l||c|c|c|c|r|}
  \hline
    \textbf{Family} &  \textbf{Median $\widetilde{\mu}$} &  \textbf{Mean $\mu$} &
    \textbf{StDev $\sigma_\mu$} &  \textbf{log(conductor)} &
    \textbf{Number}\\
    \hline \hline
\ \ 1: [0,1,1,1,T]   & 1.28  &  1.33  &  0.26 & [3.93, 9.66] &  7 \\
\ \ 2: [1,0,0,1,T]   & 1.39  &  1.40  &  0.29 & [4.66, 9.94] & 11 \\
\ \ 3: [1,0,0,2,T]   & 1.40  &  1.41  &  0.33 & [5.37, 9.97] & 11 \\
\ \ 4: [1,0,0,-1,T]  & 1.50  &  1.42  &  0.37 & [4.70, 9.98] & 20\\
\ \ 5: [1,0,0,-2,T]  & 1.40  &  1.48  &  0.32 & [4.95, 9.85] & 11\\
\ \ 6: [1,0,0,T,0]   & 1.35  &  1.37  &  0.30 & [4.74, 9.97] & 44\\
\ \ 7: [1,0,1,-2,T]  & 1.25  &  1.34  &  0.42 & [4.04, 9.46] & 10\\
\ \ 8: [1,0,2,1,T]   & 1.40  &  1.41  &  0.33 & [5.37, 9.97] & 11\\
\ \ 9: [1,0,-1,1,T]  & 1.39  &  1.32  &  0.25 & [7.45, 9.96] & 9\\
10: [1,0,-2,1,T]  & 1.34 &   1.34 &   0.42 & [3.26, 9.56] & 9\\
11: [1,1,-2,1,T]  & 1.21 &   1.19 &   0.41 & [5.73, 9.92] & 6\\
12: [1,1,-3,1,T]  & 1.32 &   1.32 &   0.32 & [5.04, 9.98] & 11\\
13: [1,-2,0,T,0]  & 1.31 &   1.29 &   0.37 & [4.73, 9.91] & 39\\
14: [-1,1,-3,1,T] & 1.45 &   1.45 &   0.31 & [5.76, 9.92] & 10\\
\hline\hline
    \textbf{All Curves}  & 1.35 &   1.36 &   0.33  &  [3.26, 9.98] &
    209\\
    \textbf{Distinct Curves} & 1.35 & 1.36 & 0.33 & [3.26, 9.98] &
    196
    \\
  \hline
\end{tabular}

\ \\ \ \\ \caption{First normalized zero above the central point for
$14$ one-parameter families of elliptic curves of rank $0$ over $\Q$
(larger conductors)}
\begin{tabular}{|l||c|c|c|c|r|}
  \hline
    \textbf{Family} &  \textbf{Median $\widetilde{\mu}$} &  \textbf{Mean $\mu$} &
    \textbf{StDev $\sigma_\mu$} &  \textbf{log(conductor)} &
    \textbf{Number}\\
    \hline \hline
\ \ 1: [0,1,1,1,T] & 0.80   & 0.86  &  0.23 &   [15.02,   15.97] &49\\
\ \ 2: [1,0,0,1,T] & 0.91   & 0.93  &  0.29 &   [15.00,   15.99] &58\\
\ \ 3: [1,0,0,2,T] & 0.90   & 0.94  &  0.30 &   [15.00,   16.00] &55\\
\ \ 4: [1,0,0,-1,T] & 0.80   & 0.90  &  0.29 &   [15.02,   16.00] &59\\
\ \ 5: [1,0,0,-2,T] &   0.75 &   0.77&    0.25 &   [15.04,   15.98]&53\\
\ \ 6: [1,0,0,T,0] & 0.75   & 0.82   & 0.27   & [15.00,   16.00]&130\\
\ \ 7: [1,0,1,-2,T]  &  0.84&    0.84 &   0.25 &   [15.04,   15.99]&63\\
\ \ 8: [1,0,2,1,T] &0.90   & 0.94 &   0.30   & [15.00,   16.00]&55\\
\ \ 9: [1,0,-1,1,T] &   0.86 &   0.89&    0.27 &   [15.02, 15.98]&57\\
10: [1,0,-2,1,T]   & 0.86 &   0.91  &  0.30  &  [15.03,   15.97]&59\\
11: [1,1,-2,1,T]   & 0.73 &   0.79  &  0.27  &  [15.00,   16.00]&124\\
12: [1,1,-3,1,T]   & 0.98 &   0.99  &  0.36  &  [15.01,   16.00]&66\\
13: [1,-2,0,T,0]   & 0.72 &   0.76  &  0.27  &  [15.00,   16.00]&120\\
14: [-1,1,-3,1,T]  & 0.90 &   0.91  &  0.24  &  [15.00,   15.99]&48\\
\hline\hline
    \textbf{All Curves}  & 0.81&  0.86 &   0.29 &   [15.00,16.00] &
    996\\
 \textbf{Distinct Curves} & 0.81 & 0.86 & 0.28 &
  [15.00,16.00] & 863\\ \hline
\end{tabular}
\end{center}
\end{table}

\newpage \ \\

\newpage \ \\

\newpage

\subsubsection{Rank 2 Curves}\label{sec:rank2inrank0fam}

The previous results were for well-separated ranges of conductors.
As the first set often has very small log-conductors, it is
possible those values are anomalous. We therefore study two sets
of curves where the log-conductors, while different, are close in
value. The goal is to see if we can detect the effect of slight
differences in the log-conductors on the repulsions.

Table $3$ provides the data from an analysis of $21$ rank $0$
one-parameter families of elliptic curves over $\Q$. The families
are from \cite{Fe2}. In each family $t$ ranges from $-1000$ to
$1000$. We searched for rank $2$ curves with log-conductor in
$[15,16]$. While we study rank $2$ curves from families of rank $2$
over $\Q$ in \S\ref{sec:famrank2oneparamoverQT2}, there the
conductors are so large that we can only analyze a few curves in
each family. In particular, there are not enough curves in one
family with conductors approximately equal to detect how slight
differences in the log-conductors effect the repulsions.

\begin{center}
\begin{table}[h]
\caption{First normalized zero above the central point for rank $2$
curves from one-parameter families of rank $0$ over $\Q$. The first
set are curves with $\log({\rm cond}) \in [15,15.5)$; the second set
are curves with $\log({\rm cond}) \in [15.5,16]$. Median =
$\widetilde{\mu}$, Mean = $\mu$, Standard Deviation (about the Mean)
= $\sigma_\mu$.}
\begin{tabular}{|l||c|c|c|r||c|c|c|r|r|}
  \hline
    \textbf{Family} &  \textbf{$\widetilde{\mu}$} &  \textbf{$\mu$} &
    \textbf{$\sigma_\mu$} &   \textbf{Number} &  \textbf{$\widetilde{\mu}$} &  \textbf{$\mu$} &
    \textbf{$\sigma_\mu$} &   \textbf{Number} \\
    \hline  \hline
\ \ 1: [0,1,3,1,T] &   1.59    &   1.83    &   0.49    &   8   &   1.71    &   1.81    &   0.40    &   19  \\
\ \ 2: [1,0,0,1,T] &   1.84    &   1.99    &   0.44    &   11  &   1.81    &   1.83    &   0.43    &   14  \\
\ \ 3: [1,0,0,2,T] &   2.05    &   2.03    &   0.26    &   16  &   2.08    &   1.94    &   0.48    &   19  \\
\ \ 4: [1,0,0,-1,T]&   2.02    &   1.98    &   0.47    &   13  &   1.87    &   1.94    &   0.32    &   10  \\
\ \ 5: [1,0,0,T,0] &   2.05    &   2.02    &   0.31    &   23  &   1.85    &   1.99    &   0.46    &   23  \\
\ \ 6: [1,0,1,1,T] &   1.74    &   1.85    &   0.37    &   15  &   1.69    &   1.77    &   0.38    &   23  \\
\ \ 7: [1,0,1,2,T] &   1.92    &   1.95    &   0.37    &   16  &   1.82    &   1.81    &   0.33    &   14  \\
\ \ 8: [1,0,1,-1,T]&   1.86    &   1.88    &   0.34    &   15  &   1.79    &   1.87    &   0.39    &   22  \\
\ \ 9: [1,0,1,-2,T]&   1.74    &   1.74    &   0.43    &   14  &   1.82    &   1.90    &   0.40    &   14  \\
10: [1,0,-1,1,T]&   2.00    &   2.00    &   0.32    &   22  &   1.81    &   1.94    &   0.42    &   18  \\
11: [1,0,-2,1,T]&   1.97    &   1.99    &   0.39    &   14  &   2.17    &   2.14    &   0.40    &   18  \\
12: [1,0,-3,1,T]&   1.86    &   1.88    &   0.34    &   15  &   1.79    &   1.87    &   0.39    &   22  \\
13: [1,1,0,T,0] &   1.89    &   1.88    &   0.31    &   20  &   1.82    &   1.88    &   0.39    &   26  \\
14: [1,1,1,1,T] &   2.31    &   2.21    &   0.41    &   16  &   1.75    &   1.86    &   0.44    &   15  \\
15: [1,1,-1,1,T]&   2.02    &   2.01    &   0.30    &   11  &   1.87    &   1.91    &   0.32    &   19  \\
16: [1,1,-2,1,T]&   1.95    &   1.91    &   0.33    &   26  &   1.98    &   1.97    &   0.26    &   18  \\
17: [1,1,-3,1,T]&   1.79    &   1.78    &   0.25    &   13  &   2.00    &   2.06    &   0.44    &   16  \\
18: [1,-2,0,T,0]&   1.97    &   2.05    &   0.33    &   24  &   1.91    &   1.92    &   0.44    &   24  \\
19: [-1,1,0,1,T]&   2.11    &   2.12    &   0.40    &   21  &   1.71    &   1.88    &   0.43    &   17  \\
20: [-1,1,-2,1,T]&   1.86    &   1.92    &   0.28    &   23  &   1.95    &   1.90    &   0.36    &   18  \\
21: [-1,1,-3,1,T]&   2.07    &   2.12    &   0.57    &   14  &   1.81    &   1.81    &   0.41    &   19  \\
\hline\hline
    \textbf{All Curves} & 1.95   &   1.97    &   0.37    &   350 &   1.85    &   1.90    &   0.40    &   388 \\
 \textbf{Distinct Curves} & 1.95 & 1.97 & 0.37 & 335 & 1.85 & 1.91
 & 0.40 & 366 \\
  \hline
\end{tabular} \normalsize
\end{table}
\end{center}

We split these rank $2$ curves from the $21$ one-parameter
families of rank $0$ over $\Q$ into two sets, those curves with
log-conductor in $[15,15.5)$ and in $[15.5,16]$. We compared the
two sets to see if we could detect the decrease in repulsion for
such small changes of the conductor. We have $21$ families, with
$350$ curves in the small conductor set and $388$ in the large
conductor set.

\begin{rek} The families are not independent: there are $15$ curves that
occur twice in the small conductor set, and $22$ in the larger. In
our amalgamations of the families we consider both the case where
we do not remove these curves, as well as the case where we do.
There is no significant difference in the results (the only
noticeable change in the table is for the mean for the larger
conductors, which increases from $1.9034$ to $1.9052$ and thus is
rounded differently). See also Remark \ref{rek:famnoindep}.
\end{rek}

The medians and means of the small conductor set are greater than
those from the large conductor set. For all curves the Pooled and
Unpooled Two-Sample $t$-Procedures give $t$-statistics of $2.5$
with over $600$ degrees of freedom; for distinct curves the Pooled
$t$-statistic is $2.16$ (respectively, the Unpooled $t$-statistic
is $2.17$) with over $600$ degrees of freedom. As the degrees of
freedom is so large, we may use the Central Limit Theorem. As
there is about a $3\%$ chance of observing a $z$-statistic of
$2.16$ or greater, the results provide evidence against the null
hypothesis (that the means are equal) at the $.05$ confidence
level, though not at the $.01$ confidence level.

While the data suggests the repulsion decreases with increasing
conductor, it is not as clear as our earlier investigations (where
we had $z$-values greater than $10$). This is, of course, due to the
closeness of the two ranges of conductors. We apply non-parametric
tests to further support our claim that the repulsion decreases with
increasing conductors. For each family in Table $3$, write a plus
sign if the small conductor set has a greater mean and a minus sign
if not. There are four minus signs and seventeen plus signs. The
null hypothesis is that each mean is equally likely to be larger.
Under the null hypothesis, the number of minus signs is a random
variable from a binomial distribution with $N=21$ and $\theta =
\foh$. The probability of observing four or fewer minus signs is
about $3.6\%$, supporting the claim of decreasing repulsion with
increasing conductor. For the medians there are seven minus signs
out of twenty-one; the probability of seven or fewer minus signs is
about $9.4\%$. Every time the smaller conductor set had the lesser
mean, it also had the lesser median; the mean and median tests are
not independent.


\newpage
\subsection{One-Parameter Families of Rank $2$ Over
$\Q$}\label{sec:famrank2oneparamoverQT2}
\subsubsection{Family $y^2 = x^3 -T^2x
+T^2$}\label{sec:exp2fam}

We study the first normalized zero above the central point for $69$
rank $2$ elliptic curves from the one-parameter family $y^2 = x^3
-T^2x+T^2$ of rank $2$ over $\Q$. There are $35$ curves with
$\log({\rm cond}) \in [7.8,16.1]$ in Figure $9$ and 34 with
$\log({\rm cond}) \in [16.2,23.3]$ in Figure $10$. Unlike the
previous examples where we chose many curves of the same rank from
different families, here we have just one family. As the conductors
grow rapidly, we have far fewer data points, and the range of the
log-conductors is much greater. However, even for such a small
sample, the repulsion decreases with increasing conductors, and the
shape begins to approach the conjectured distribution. The Pooled
and Unpooled Two-Sample $t$-Procedures give $t$-statistics over $5$
with over $60$ degrees of freedom, and we may use the Central Limit
Theorem. As the probability of a $z$-value of $5$ or more is less
than $10^{-4}$ percent, the data does not support the null
hypothesis (i.e., the data supports our conjecture that the
repulsion decreases as the conductors increase).

\begin{figure}[h]
\begin{center}
\includegraphics[width=7cm]{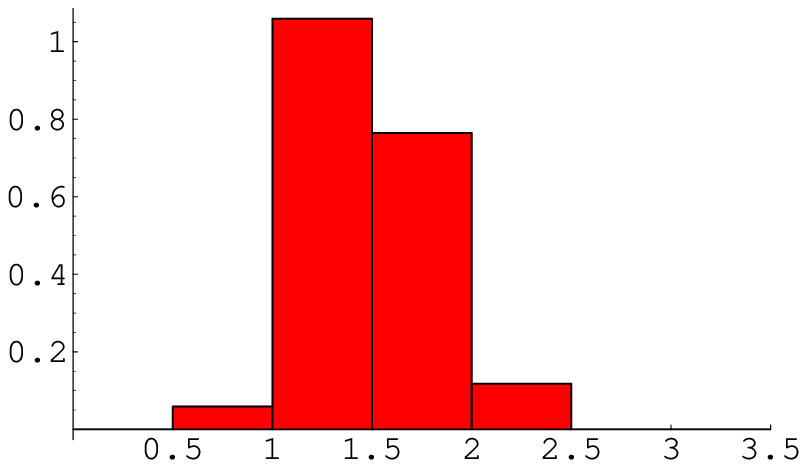}\\
Figure $9$: First normalized zero above the central point\\ from
rank $2$ curves in the family $y^2 = x^3 - T^2 x + T^2$.\\ $35$
curves, $\log({\rm cond}) \in [7.8, 16.1]$, $\text{median} =
1.85$, $\text{mean} = 1.92$,\\ standard deviation about the mean
$=.41$
\end{center}
\end{figure}

\begin{figure}[h]
\begin{center}
\includegraphics[width=7cm]{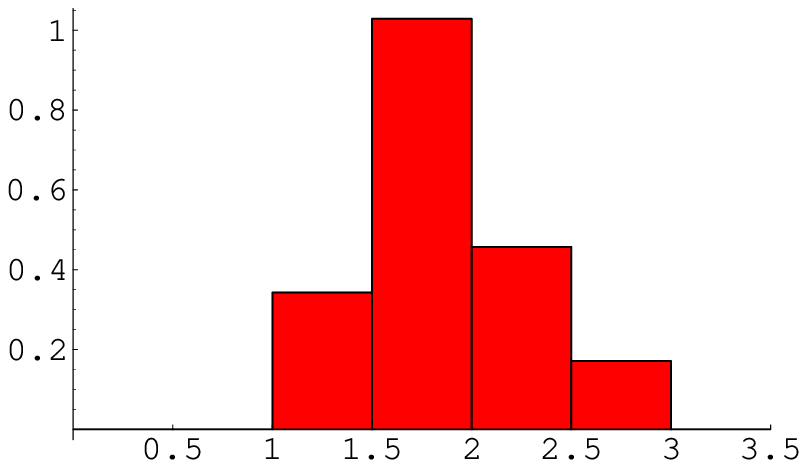}\\
Figure $10$: First normalized zero above the central point\\ from
rank $2$ curves in the family $y^2 = x^3 - T^2 x + T^2$.\\ $34$
curves, $\log({\rm cond})\in [16.2, 23.3]$, $\text{median} =
1.37$, $\text{mean} = 1.47$,\\ standard deviation about the mean
$=.34$
\end{center}
\end{figure}

\newpage
\subsubsection{Rank 2 Curves}

We consider $21$ one-parameter families of rank $2$ over $\Q$, and
investigate curves of rank $2$ in these families. The families are
from \cite{Fe2}. We again amalgamated the different families, and
summarize the results in Table 4.

The difference between these experiments and those of
\S\ref{sec:rank2inrank0fam} is that, while both deal with
one-parameter families over $\Q$, here we study curves of rank $2$
from families of rank $2$ over $\Q$; earlier we studied curves of
rank $2$ from families of rank $0$ over $\Q$. If the Density
Conjecture (with orthogonal groups) is correct for the entire
one-parameter family, in the limit 0\% of the curves in a family of
rank $r$ have rank $r+2$ or greater. Thus our previous
investigations of curves of rank $2$ in a family of rank $0$ over
$\Q$ were a study of a measure zero subset. Unlike curves of rank
$2$ in families of rank $2$ over $\Q$, we have no theoretical
evidence supporting a proposed random matrix model for curves of
rank $2$ in families of rank $0$. We compare the results from rank
$2$ curves in rank $2$ families over $\Q$ to the rank $2$ curves
from rank $0$ families over $\Q$ in \S\ref{sec:compareoneparams}.

\begin{center}
\begin{table}[h]
\caption{First normalized zero above the central point for 21
one-parameter families of rank $2$ over $\Q$ with $\log({\rm
cond}) \in [15,16]$ and $t \in [0,120]$. The median of the first
normalized zero of the 64 curves is 1.64.}
\begin{tabular}{|l||c|c|c|c|}
 \hline
\textbf{Family}      &   \textbf{Mean}    &   \textbf{Standard Deviation}
&  \textbf{log(conductor)}  &  \textbf{Number} \\
\hline
\ \ 1:   [1,T,0,-3-2T,1] &   1.91    &   0.25    &  [15.74,16.00] &\ \ 2   \\
\ \ 2:   [1,T,-19,-T-1,0]    &   1.57    &   0.36  & [15.17,15.63]  &\ \    4   \\
\ \ 3:   [1,T,2,-T-1,0]  &   1.29    &    & [15.47, 15.47] & \ \   1   \\
\ \ 4:   [1,T,-16,-T-1,0]    &   1.75    &   0.19    &  [15.07,15.86] &\ \  4   \\
\ \ 5:   [1,T,13,-T-1,0] &   1.53    &   0.25    &  [15.08,15.91] &\ \  3   \\
\ \ 6:   [1,T,-14,-T-1,0]    &   1.69    &   0.32    &  [15.06,15.22] &\ \  3   \\
\ \ 7:   [1,T,10,-T-1,0] &   1.62    &   0.28    &  [15.70,15.89] &\ \  3   \\
\ \ 8:   [0,T,11,-T-1,0] &   1.98    &    &   [15.87,15.87] &\ \  1   \\
\ \ 9:   [1,T,-11,-T-1,0]    &   &    &   &   \\
10:  [0,T,7,-T-1,0]  &   1.54    &   0.17    & [15.08,15.90] & \ \  7   \\
11:  [1,T,-8,-T-1,0] &   1.58    &   0.18    & [15.23,25.95] &\ \   6   \\
12:  [1,T,19,-T-1,0] &    &       &  &    \\
13:  [0,T,3,-T-1,0]  &   1.96    &   0.25    & [15.23, 15.66] & \ \  3   \\
14:  [0,T,19,-T-1,0] &    &    &    &  \\
15:  [1,T,17,-T-1,0] &   1.64    &   0.23    & [15.09, 15.98] & \ \  4   \\
16:  [0,T,9,-T-1,0]  &   1.59    &   0.29    & [15.01, 15.85] & \ \  5   \\
17:  [0,T,1,-T-1,0]  &   1.51    &    &   [15.99, 15.99] &\ \  1   \\
18:  [1,T,-7,-T-1,0] &   1.45    &   0.23    &   [15.14, 15.43] &\ \  4   \\
19:  [1,T,8,-T-1,0]  &   1.53    &   0.24    & [15.02, 15.89] &  10  \\
20:  [1,T,-2,-T-1,0] &   1.60    &    &  [15.98, 15.98] & \ \ 1   \\
21:  [0,T,13,-T-1,0] &   1.67    &   0.01    & [15.01, 15.92] &\ \   2   \\
  \hline \hline \textbf{All Curves} & 1.61 & 0.25 & [15.01, 16.00] & 64 \\ \hline
\end{tabular}
\end{table}
\end{center}

\begin{rek} There
are $23$ rank $4$ curves in the 21 one-parameter families of rank
$2$ over $\Q$ with log-conductors in $[15, 16]$ and $t \in [0,
120]$. For the first normalized zero above the central point, the
median is 3.03, the mean is 3.05, and the standard deviation about
the mean is 0.30.
\end{rek}


\newpage
\subsection{Comparison Between One-Parameter Families of Different
Rank}\label{sec:compareoneparams}

In Table 5 we investigate how the first normalized zero above the
central point of rank $2$ curves depends on how the curves are
obtained. The first family is rank $2$ curves from the 21
one-parameter families of rank $0$ over $\Q$ from Table 3, while the
second is rank $2$ curves from the 21 one-parameter families of rank
$2$ over $\Q$ from Table 4; in both sets the log-conductors are in
$[15,16]$. A $t$-Test on the two means gives a $t$-statistic of
6.60, indicating the two means differ. Thus the mean of the first
normalized zero above the central point of rank $2$ curves in a
one-parameter family over $\Q$ (for conductors in this range)
depends on \emph{how} we choose the curves. For the range of
conductors studied, rank $2$ curves from rank $0$ one-parameter
families over $\Q$ do \emph{not} behave the same as rank $2$ curves
from rank $2$ one-parameter families over $\Q$.

\begin{table}[h]
\begin{center}
\caption{First normalized zero above the central point. The first
family is the 701 rank $2$ curves from the $21$ one-parameter
families of rank $0$ over $\Q$ from Table 3 with $\log({\rm cond})
\in [15,16]$; the second family is the 64 rank $2$ curves from the
$21$ one-parameter families of rank $2$ over $\Q$ from Table $4$
with $\log({\rm cond}) \in [15,16]$.}
\begin{tabular}{|l||c|c|c|c|}
  \hline
\textbf{Family} & \textbf{Median} & \textbf{Mean} & \textbf{Std.
Dev.} & \textbf{Number} \\ \hline
Rank 2 Curves, Families Rank 0 over $\Q$ & 1.926 & 1.936 & 0.388 & 701 \\
Rank 2 Curves, Families Rank 2 over $\Q$ & 1.642 & 1.610 & 0.247 & \ \ 64 \\
\hline
\end{tabular}
\end{center}
\end{table}

\newpage
\subsection{Spacings between normalized zeros}\label{sec:spacingsnorzeros}

For finite conductors, even when there are no zeros at the central
point, the first normalized zero above the central point is
repelled from the predicted $N\to\infty$ scaling limits. The
repulsion increases with the number of zeros at the central point
and decreases with increasing conductor. For an elliptic curve
$E$, let $z_1, z_2, z_3, \dots$ denote the imaginary parts of the
normalized zeros above the central point. We investigate whether
or not $z_{j+1} - z_j$ depends on the repulsion from the central
point.

We consider the following two sets of curves in Table 6: \bi \item
the $863$ distinct rank $0$ curves with $\log({\rm cond}) \in
[15,16]$ from the $14$ one-parameter families of rank $0$ over $\Q$
from Table $2$; \item the $701$ distinct rank $2$ curves with
$\log({\rm cond}) \in [15,16]$ from the $21$ one-parameter families
of rank $0$ over $\Q$ from Table $3$. \ei In Table $6$ we calculate
the median and mean for $z_2-z_1$, $z_3-z_2$ and $z_3-z_1$. The last
statistic involves the sum of differences between two adjacent
normalized zeros, and allows the possibility of some effects being
averaged out. While the normalized zeros are repelled from the
central point (and by different amounts for the two sets), the
\emph{differences} between the normalized zeros are statistically
independent of this repulsion. We performed a $t$-Test on the means
in the three cases. For each case the $t$-statistic was less than
$2$, strongly supporting the null hypothesis that the differences
are independent of the repulsion.

\begin{center}
\begin{table}[h]
\caption{Spacings between normalized zeros. All curves have
$\log({\rm cond}) \in [15,16]$, and $z_j$ is the imaginary part of
the $j$\textsuperscript{th} normalized zero above the central
point. The $863$ rank 0 curves are from the $14$ one-parameter
families of rank $0$ over $\Q$ from Table $2$; the $701$ rank $2$
curves are from the $21$ one-parameter families of rank $0$ over
$\Q$ from Table $3$.}
\begin{tabular}{|l||c|c||r|}
  \hline
& \textbf{863 Rank $0$ Curves} & \textbf{701 Rank $2$ Curves}
& \textbf{t-Statistic} \\
\hline
\textbf{Median $z_2-z_1$}         & 1.28 & 1.30 & \ \\
\textbf{Mean\ \ \ \ $z_2 - z_1$}  & 1.30 & 1.34 & -1.60 \\
\textbf{StDev\ \ \ \ $z_2 - z_1$} & 0.49 & 0.51 & \ \\
\hline
\textbf{Median $z_3-z_2$}          & 1.22 & 1.19 & \ \\
\textbf{Mean\ \ \ \  $z_3 - z_2$}  & 1.24 & 1.22 & 0.80 \\
\textbf{StDev\ \ \ \  $z_3 - z_2$} & 0.52 & 0.47 & \ \\
\hline
\textbf{Median $z_3-z_1$}          & 2.54 & 2.56 & \ \\
\textbf{Mean\ \ \ \  $z_3 - z_1$}  & 2.55 & 2.56 & -0.38\\
\textbf{StDev\ \ \ \  $z_3 - z_1$} & 0.52 & 0.52 & \ \\

\hline
\end{tabular}
\end{table}
\end{center}

We have consistently observed that the more zeros at the central
point, the greater the repulsion. One possible explanation is as
follows: for rank $2$ curves in a rank $0$ one-parameter family
over $\Q$, the first zero above the central point collapses down
to the central point, and the other zeros are left alone. As the
zeros are symmetric about the central point, the effect of one
zero above the central point collapsing is to increase the number
of zeros at the central point by $2$.

For our $14$ one-parameter families of elliptic curves of rank $0$
over $\Q$ and log-conductors in $[15,16]$, we studied the second and
third normalized zero above the central point. The mean of the
second normalized zero is $2.16$ with a standard deviation of $.39$,
while the third normalized zero has a mean of $3.41$ and a standard
deviation of $.41$. These numbers statistically differ\footnote{The
Pooled and Unpooled $t$-statistics in both experiments are greater
than $6$, providing evidence against the null hypothesis that the
two means are equal.} from the first and second normalized zeros of
the rank $2$ curves from our $21$ one-parameter families of rank $0$
over $\Q$ with log-conductor in $[15,16]$, where the means were
respectively $1.93$ (with a standard deviation of $.39$) and $3.27$
(with a standard deviation of $.39$). Thus while for a given range
of log-conductors the average second normalized zero of a rank $0$
curve is close to the average first normalized zero of a rank $2$
curve, they are not equal and the additional repulsion from extra
zeros at the central point cannot be entirely explained by
\emph{only} collapsing the first zero to the central point while
leaving the other zeros alone.

\begin{rek}\label{rek:attraction}
As the second (resp., third) normalized zero for rank 0 curves
in rank 0 families over $\Q$ is 2.16 (resp., 3.41) while the first
(resp., second) normalized zero for rank 2 curves in rank 0 families
over $\Q$ is 1.93 (resp., 3.27), one can interpret the effect of the
additional zeros at the central point as an \emph{attraction}.
Specifically, for curves of rank 2 in a rank 0 family over $\Q$, by
symmetry two zeros collapse to the central point, and the remaining
zeros are then attracted to the central point, being closer than the
corresponding zeros from rank 0 curves. As remarked in \S3.5 of
\cite{Far}, the term ``lowest zero'' is not well defined when there
are multiple zeros at the central point. We can either mean the
first zero above the central point, or one of the many zeros at the
central point. In all cases, for finite conductors there is
repulsion from the $N\to\infty$ scaling limits of random matrix
theory; however, ``attraction'' might be a better term for the
effect of additional zeros at the central point, though the current
terminology is to talk about repulsion of zeros at the central
point.
\end{rek}

We now study the differences between normalized zeros coming from
one-parameter families of rank $2$ over $\Q$. Table $7$ shows that
while the normalized zeros are repelled from the central point, the
\emph{differences} between the normalized zeros are statistically
independent of the repulsion. We performed a $t$-Test for the means
in the three cases studied. For two of the three cases the
$t$-statistic was less than $2$ (and in the third it was only
$2.05$), supporting the null hypothesis that the differences are
independent of the repulsion.

\begin{table}[h]
\begin{center}
\caption{Spacings between normalized zeros. All curves have
$\log({\rm cond}) \in [15,16]$, and $z_j$ is the imaginary part of
the $j$\textsuperscript{th} normalized zero above the central point.
The $64$ rank $2$ curves are the $21$ one-parameter families of rank
$2$ over $\Q$ from Table $4$; the $23$ rank $4$ curves are the $21$
one-parameter families of rank $2$ over $\Q$ from Table $4$.}
\begin{tabular}{|l||c|c||r|}
  \hline
& \textbf{64 Rank $2$ Curves} & \textbf{23 Rank $4$ Curves}
& \textbf{t-Statistic} \\
\hline
\textbf{Median $z_2-z_1$}         & 1.26 & 1.27 &  \\
\textbf{Mean\ \ \ \ $z_2 - z_1$}  & 1.36 & 1.29 & 0.59 \\
\textbf{StDev\ \ \ \ $z_2 - z_1$} & 0.50 & 0.42 &  \\
\hline
\textbf{Median $z_3-z_2$}          & 1.22 & 1.08 &  \\
\textbf{Mean\ \ \ \  $z_3 - z_2$}  & 1.29 & 1.14 & 1.35 \\
\textbf{StDev\ \ \ \  $z_3 - z_2$} & 0.49 & 0.35 &  \\
\hline
\textbf{Median $z_3-z_1$}          & 2.66 & 2.46 &  \\
\textbf{Mean\ \ \ \  $z_3 - z_1$}  & 2.65 & 2.43 & 2.05\\
\textbf{StDev\ \ \ \  $z_3 - z_1$} & 0.44 & 0.42 &  \\
\hline
\end{tabular}
\end{center}
\end{table}

We performed one last experiment on the differences between
normalized zeros. In Table $8$ we compare two sets of rank $2$
curves: the first are the $21$ one-parameter families of rank $0$
over $\Q$ from Table $3$, while the second are the $21$
one-parameter families of rank $2$ over $\Q$ from Table $4$. While
the first normalized zero is repelled differently in the two
cases, the differences are statistically independent from the
nature of the zeros at the central point, as indicated by all
$t$-statistics being less than $2$. This suggests that the
\emph{spacings} between adjacent normalized zeros above the
central point is independent of the repulsion at the central
point; in particular, this quantity does not depend on how we
construct our family of rank $2$ curves.

\begin{table}
\begin{center}
\caption{Spacings between normalized zeros. All curves have
$\log({\rm cond}) \in [15,16]$, and $z_j$ is the imaginary part of
the $j$\textsuperscript{th} normalized zero above the central
point. The 701 rank 2 curves are the $21$ one-parameter families
of rank $0$ over $\Q$ from Table $3$, and the 64 rank $2$ curves
are the $21$ one-parameter families of rank $2$ over $\Q$ from
Table $4$.}
\begin{tabular}{|l||c|c||r|}
  \hline
& \textbf{701 Rank $2$ Curves} & \textbf{64 Rank $2$ Curves}
& \textbf{t-Statistic} \\
\hline
\textbf{Median $z_2-z_1$}         & 1.30 & 1.26 & \\
\textbf{Mean\ \ \ \ $z_2 - z_1$}  & 1.34 & 1.36 & 0.69\\
\textbf{StDev\ \ \ \ $z_2 - z_1$} & 0.51 & 0.50 & \\
\hline
\textbf{Median $z_3-z_2$}          & 1.19 & 1.22 & \\
\textbf{Mean\ \ \ \  $z_3 - z_2$}  & 1.22 & 1.29 & 1.39 \\
\textbf{StDev\ \ \ \  $z_3 - z_2$} & 0.47 & 0.49 & \\
\hline
\textbf{Median $z_3-z_1$}          & 2.56 & 2.66 & \\
\textbf{Mean\ \ \ \  $z_3 - z_1$}  & 2.56 & 2.65 & 1.93\\
\textbf{StDev\ \ \ \  $z_3 - z_1$} & 0.52 & 0.44 & \\

\hline
\end{tabular}
\end{center}
\end{table}


\section{Summary and future work}

As the conductors tend to infinity, theoretical results support the
validity of the $N\to\infty$ scaling limit of the Independent Model
for all curves in one-parameter families of elliptic curves of rank
$r$ over $\Q$; however, it is unknown what the correct model is for
the sub-family of curves of rank $r+2$. The experimental evidence
suggests that the first normalized zero, for small and finite
conductors, is repelled by zeros at the central point. Further, the
more zeros at the central point, the greater the repulsion; however,
the repulsion decreases as the conductors increase, and the
difference between adjacent normalized zeros is statistically
independent of the repulsion and the rank of the curves.

At present we can calculate the first normalized zero for
log-conductors about $25$. While we can use more powerful
computers to study larger conductors, it is unlikely these
conductors will be large enough to see the predicted limiting
behavior. It is interesting that, unlike the excess rank
investigations, we see noticeable convergence to the limiting
theoretical results as we increase the conductors.

An interesting project is to determine a theoretical model to
explain the behavior for finite conductors. In the large-conductor
limit, analogies with the function field and calculations with the
explicit formula lead us to the Independent Model for curves of rank
$r$ from families of rank $r$ over $\Q$, and theoretical results in
the number field case support this. It is not unreasonable to posit
that in the finite-conductor analogues the size of the matrices
should be a function of the log-conductors. Unfortunately the
statistics for the finite $N\times N$ random matrix ensembles are
expressed in terms of eigenvalues of integral equations, and are
usually only plotted in the $N\to\infty$ scaling limit. This makes
comparison with the experimental data difficult, and a future
project is to analyze the finite $N$ cases by using the finite $N$
kernels. Such an analysis will facilitate comparing the finite $N$
limits of the Independent and Interaction Models for curves of rank
$r+2$ from families of rank $r$ over $\Q$.

\section*{Acknowledgement}

The program used to calculate the order of vanishing of elliptic
curve $L$-functions at the central point was written by Jon Hsu, Leo
Goldmakher, Stephen Lu and the author; the programs to calculate the
first few zeros above the central point in families of elliptic
curves is due to Adam O'Brien, Aaron Lint, Atul Pokharel, Michael
Rubinstein and the author. I would like to thank Eduardo Due\~nez,
Frank Firk, Michael Rosen, Peter Sarnak, Joe Silverman and Nina
Snaith for many enlightening conversations, the referees for a very
thorough reading of the paper and suggestions which improved the
presentation, and the Information Technology Managers at the
Mathematics Departments at Princeton, The Ohio State University and
Brown for help in getting all the programs to run compatibly.

\newpage

\appendix


\section{``Harder'' Ensembles of Orthogonal Matrices}\label{app:eduardo}\label{sec:hard-edge-ensembl}

\centerline{Eduardo Due\~nez} \centerline{Department of Applied
  Mathematics} \centerline{University of Texas at San Antonio}
\centerline{6900 N Loop 1604 W}
\centerline{San Antonio, TX 78249}
\centerline{\texttt{eduenez@math.utsa.edu}} \ \\

In this appendix we derive the conditional (interaction) eigenvalue
probability measure~\eqref{eq:14} and illustrate how it affects
eigenvalue statistics near the central point~$1$, in particular
through repulsion (observed via the $1$-level density).  We also
explain the relation to the classical Bessel kernels of random matrix
theory, and to other central-point statistics.

\subsection{Full Orthogonal Ensembles}
\label{sec:full-SO}

In view of our intended application we will be concerned exclusively
with random matrix ensembles of orthogonal matrices in what follows.
If we write the eigenvalues (in no particular order) of a special%
\footnote{That is, of determinant one.} orthogonal matrix of size~$2N$
(resp.,~$2N+1$) as $\{\pm e^{i\theta_j}\}_1^N$ (resp.,
$\{+1\}\cup\{\pm e^{i\theta_j}\}_1^N$) with $0\leq \theta_j\leq \pi$
then the $N$-tuple $\Theta = $ $(\theta_1,\dots,\theta_N)$
parametrizes the
eigenvalues. 
In terms of the angles $\theta_j$, the probability measure of the
eigenvalues induced from normalized Haar measure on SO($2N$)
(resp., on SO($2N+1$) upon discarding one forced eigenvalue of
$+1$) can be identified with a measure on $[0,\pi]^N$,
\begin{align}
  \label{eq:7b}
  d\ve_0(\Th) &= \tilde C_{N}^{(0)}
  \prod_{1\leq j<k\leq N}(\cos\theta_k-\cos\theta_j)^2\prod_{1\leq j\leq N} d\theta_j \\
  \label{eq:8}
  d\ve_1(\Th) &= \tilde C_{N}^{(1)} \prod_{1\leq j<k\leq
    N}(\cos\theta_k-\cos\theta_j)^2 \prod_{1\leq j\leq N}
  \sin^2({\textstyle\frac{\theta_j}2})d\theta_j
\end{align}
in the $2N$ and $2N+1$ cases, respectively, as shown in
\cite{Con,KaSa1}; the normalization constants $\tilde C_N^{(m)}$
ensure that the measures on the right-hand sides are probability
measures. Note that formulas~\eqref{eq:7b} and~\eqref{eq:8} are
symmetric upon permuting the $\theta_j$'s, so issues related to a
choice of a particular ordering of the eigenvalues of the matrix
are irrelevant.  More importantly, observe the quadratic exponent
of the differences of the cosines.

The statistical behavior of the eigenvalues near~$+1$ is closely
related to the order of vanishing of the measures above at
$\theta=0$. We change variables and replace the eigenvalues
$e^{\pm i\theta}$ by the levels
\begin{equation}
  \label{eq:15}
  x\ = \ \cos\theta
\end{equation}
so the measures above become measures on $[-1,+1]^N$:
\begin{align}
  \label{eq:9}
  d\ve_0(X) &\ = \ C_{N}^{(0)} \prod_{1\leq j<k\leq N}(x_k-x_j)^2\prod_{j=1}^N (1-x_j)^\mhalf(1+x_j)^\mhalf dx_j \\
  \label{eq:10}
  d\ve_1(X) &\ = \ C_{N}^{(1)} \prod_{j<k}(x_k-x_j)^2\prod_{j=1}^N
  (1-x_j)^\half(1+x_j)^\mhalf dx_j,
\end{align}
where $X=(x_1,\dots,x_N)$ and $C_N^{(m)}$ are suitable
normalization constants.  Here we observe the appearance of the
weight functions on $[-1,1]$
\begin{align}
  \label{eq:11}
  w(x) = (1-x)^a(1+x)^\mhalf,\qquad a=
  \begin{cases}
    -1/2 & \text{for SO($2N$)}\\
    +1/2 & \text{for SO($2N+1$).}
  \end{cases}
\end{align}
By the Gaudin-Mehta theory (see for example \cite{Meh}), and in view
of the quadratic exponent of the differences of the ``levels''
$x_j$, the study of eigenvalue statistics using classical methods is
intimately related to the sequence of orthogonal polynomials with
respect to the
weight $w(x)$.%
\footnote{The inner product being $\langle f,g\rangle = \int_0^1
  f(x)\overline{g(x)}w(x)dx$.}


In classical random matrix theory terminology (especially in the
context of the Laguerre and Jacobi ensembles) the endpoints
$-1,+1$ are called the ``hard edges'' of the spectrum because the
probability measure, considered on $\R^N$, vanishes outside
$[-1,+1]^N$.  We will keep calling $\theta=0$, $\pi$ the ``central
points'' (endpoints of the diameter with respect to which the
spectrum is symmetric). Phenomena about central points and hard
edges are equivalent in view of the change of
variables~\eqref{eq:15}.  Perhaps not surprisingly, the
parameter~$a$, which dictates the order of vanishing of the weight
function $w(x)$ at the hard edge~$+1$, suffices to characterize
the mutually different statistics near the central point in each
of SO(even) and~SO(odd).  However, the importance of this
parameter is best understood in the context of certain
sub-ensembles of SO as described below.

\subsection{Conditional (``Harder'') Orthogonal Ensembles}
\label{sec:Hard-SO}

The conditional eigenvalue measure for the sub-ensemble
$SO^{(2r)}(2N)$ of $SO(2N)$ consisting of matrices for which some
$2r$ of the $2N$~eigenvalues are equal to~$+1$ can easily be
obtained from~\eqref{eq:9}.  Let
\begin{equation}
  f(x_1,\dots,x_N) =
  C_N^{(m)}\prod_{1\leq j< k\leq N}(x_k-x_j)^2\prod_{1\leq j\leq N} w(x_j)
\end{equation}
be the normalized probability density function of the levels for
$SO(2N)$, where $w(x)$ is as in~\eqref{eq:11} with $a=-1/2$ and
$m=0$.
Now let $t_1,\dots,t_r$ be chosen so $0<t_k<1$, let
$K=\prod_j[1-t_j,1]$ and $I=J\times K$ for some box
$J\subset[-1,1]^{N-r}$.  This means we are constraining $r$ pairs
of levels to lie in a neighborhood of $x=1$ (or equivalently that
we are construing $r$ pairs of eigenvalues to lie in circular
sectors about the point $1$ on the unit circle).  Thus, the
conditional probability that the remaining $N-r$ pairs of
eigenvalues lie in~$J$ is given by
\begin{equation}
  \label{eq:condl}
  F(T;J)\ =\ \frac{\displaystyle\int_{J\times K}f(x)dx}
{\displaystyle\int_{[-1,1]^{N-r}\times K}f(x)dx},
\end{equation}
where $T=(t_1,\dots,t_r)$.  The conditional probability measure of
the eigenvalues for the sub-ensemble $SO^{(2r)}(2N)$ is the limit
as all $t_k\to0+$ of $F(t;J)$, as a function of the box $J$, call
it $G(J)$.  Applying L'H\^opital's rule $r$ times to the
quotient~\eqref{eq:condl} (once on each variable $t_k$) and using
the fundamental theorem of calculus we get
\begin{equation}
  \label{eq:G}
  G(J) = \lim_{T\to0+}
  \frac{\int_J (\Van(X))^2 (M(X,T))^2 w(X)
    dX\cdot (\Van(T))^2w(T)}
  {\int_{[0,1]^{N-r}} (\Van(X))^2 (M(X,T))^2 w(X)
    dX\cdot (\Van(T))^2w(T)}
\end{equation}
where $X=(x_1,\dots,x_{N-r})$, and
\begin{align*}
  \Van(X)&=\prod_{1\leq j<k\leq N-r}(x_k-x_j) &
  \Van(T)&=\prod_{1\leq j<k\leq r}(t_k-t_j)  \\
  w(X)&=\prod_{1\leq j\leq N-r}w(x_j) &
  w(T)&=\prod_{1\leq k\leq r}w(1-t_k) \\
  M(X,T)&=\prod_{1\leq j\leq N-r\atop 1\leq k\leq r}(1-t_k-x_j).
\end{align*}
Naturally, the factors of $\Van(T)$ and $w(T)$ cancel in
equation~\eqref{eq:G}.  Since $M(X,T)$ is bounded, the integrands
in equation~\eqref{eq:G} are uniformly dominated by an integrable
function, ensuring that we can let all $t_j\to0$ in the integrands
of~\eqref{eq:G} to obtain
\begin{equation}
  \label{eq:GG}
  G(J) =  \frac{\int_J (\Van(X))^2 (M(X,0))^2 w(X)
    dX}
  {\int_{[0,1]^{N-r}} (\Van(X))^2 (M(X,0))^2 w(X)
    dX}.
\end{equation}
Now observe that $(M(X,0))^2w(X)=\prod_{1\leq j\leq r}\widetilde
w(x_j)$ where $\widetilde w(x)$ is given by equation~\eqref{eq:11}
with $a$ replaced by $\tilde a=a+2$, so the probability measure of
the eigenvalues for $SO^{(2r)}(2N)$ is obtained from that of
$SO(2(N-m))$ simply by changing the weight function
$w\mapsto\widetilde w$. Explicitly, the probability measure of the
eigenvalues of the ensemble $SO^{(2r)}(2N)$ is given by
\begin{equation}
  \label{eq:14b}
  d\varepsilon_{m}(X) = C_{N-m}^{(m)}
  \prod_{j<k}(x_k-x_j)^2\prod_j(1-x_j)^{m-\frac12}
  \prod_j dx_j,
\end{equation}
where $m=2r$, $X=(x_1,\dots,x_{N-r})$, the indices $j,k$ range
from $1$ to $N-r$, and $C_{N-m}^{(m)}$ are suitable normalization
constants (equal to the reciprocal of the denominator of the
right-hand side of equation~\eqref{eq:GG}.)

The same argument shows that the sub-ensemble $SO^{(2r+1)}(2N+1)$
of $SO(2N+1)$ consisting of matrices for which $2r+1$ eigenvalues
are equal to~$+1$ has the same eigenvalue measure~\eqref{eq:14b}
with $m=2r+1$.  Because the density of the measure vanishes to a
higher order near the edge~$+1$ the larger $m$ is, we will say
that the edge becomes harder when $m$ is larger (whence the title
of this section), and call $m$ its hardness.


\subsection{Independent Model}
\label{sec:indep-mod}

It is important to observe that the presence of the $m$-multiple
eigenvalues at the central point in these harder sub-ensembles of
orthogonal matrices has a strong repelling effect due to the extra
factor $(1-x)^m$ multiplied by the weight
$(1-x)^\mhalf(1+x)^\mhalf$ of SO(even).  For comparison purposes
consider the following situation, first in the SO(even) case.  The
number of eigenvalues equal to~$+1$ of any $SO(2N)$ matrix is
always an even number~$2r$, and one may consider the sub-ensemble
$\mathcal{A}_{2N,2r}$ of $SO(2N)$,
\begin{equation}
  \label{eq:28b}
  \mathcal{A}_{2N,2r} = \left\{
    \begin{pmatrix}
      I_{2r \times 2r} & \\
      & g
    \end{pmatrix}
    : g \in SO(2N-2r)\right\},
\end{equation}
which is just $SO(2N-2r)$ in disguise.  This is certainly a
sub-ensemble of $SO(2N)$ consisting of matrices with at least $2r$
eigenvalues equal to~$+1$, albeit quite a different one from the
$2r$-hard sub-ensemble of $SO(2N)$ described before.  For example,
the eigenvalue measure (apart from the point masses at the last
$2r$ eigenvalues) for $\mathcal{A}_{2N,2r}$ is
\begin{equation}
  \label{eq:29}
  d\varepsilon_0(x_1,\dots,x_{N-r}),
\end{equation}
and not $d\varepsilon_{2r}(x_1,\dots,x_{N-r})$.  The same
observation applies in the $SO(2N+1)$ case: with the obvious
notation, the
eigenvalue measure for $\mathcal{A}_{2N+1,2r}$ is%
\footnote{Observe that a matrix in $\mathcal{A}_{2N+1,2r}$ has
$2r+1$
  eigenvalues equal to $+1$ and $N-r$ other pairs eigenvalues.}
\begin{equation}
  \label{eq:30}
  d\varepsilon_1(\theta_1,\dots,\theta_{N-r}),
\end{equation}
and not $d\varepsilon_{2r+1}(\theta_1,\dots,\theta_{N-r})$.

\subsection{$1$-Level Density: Full Orthogonal}
\label{sec:1-level-full}

Before dealing with the harder sub-ensembles, we make some
comments about the hard edges of the full SO(even) and SO(odd).
The local statistics near the point $+1$ are dictated by the even
`$+$' (resp., odd `$-$') Sine Kernels
\begin{equation}
  S_{\pm}(\xi,\eta) = S(\xi,\eta)\pm S(\xi,-\eta)
\end{equation}
in the case of SO(even) (resp., SO(odd)); see \cite{KaSa1,KaSa2}.
Here $\xi,\eta$ are rescaled variables centered about the value $0$,
namely related to the original variables by%
\footnote{This is justified by the fact $N/\pi$ is the average
  (angular) \emph{asymptotic} density of the eigen-angles $\theta_j$
  of a random orthogonal matrix, hence asymptotic equidistribution holds
  ---away from the central points!}
\begin{equation}
  \label{eq:xi}
  x = \cos\left(\frac\pi N\xi\right)
\end{equation}
and $S(x,y)=\sin(\pi x)/(\pi x)$ is the Sine Kernel, which has the
universal property of describing the local statistics at
\emph{any} bulk point of \emph{any} ensemble with local quadratic
local level repulsion~\cite{DKMcVZ}. However, it is not the Sine
Kernel but its even (resp., odd) counterparts that dictate the
local statistics near the central point.
For example, the central one-level density is given by the
diagonal values at $x=y$ of the respective kernel:
\begin{align}
  \label{eq:4}
  \rho_+(x) &= 1 + \frac{\sin2\pi x}{2\pi x},& \qquad\text{for SO(even),} \\
  \label{eq:5}
  \rho_-(x) &= 1 + \frac{\sin2\pi x}{2\pi x} + \delta(x),&
  \qquad\text{for SO(odd).}
\end{align}
(In the SO(odd) case the Dirac delta reflects the fact that any
such matrix has an eigenvalue at the central point.)  Observe that
$\rho_-$ vanishes to second order, whereas $\rho_+$ does not
vanish at the
central point $x=0$.%
\footnote{If the central point were not atypical, then the local
  density would be dictated by the diagonal values $S(x,x)\equiv1$ of
  the Sine Kernel.}
\begin{figure}[!htb]
  \centering
  \includegraphics[scale=.9]{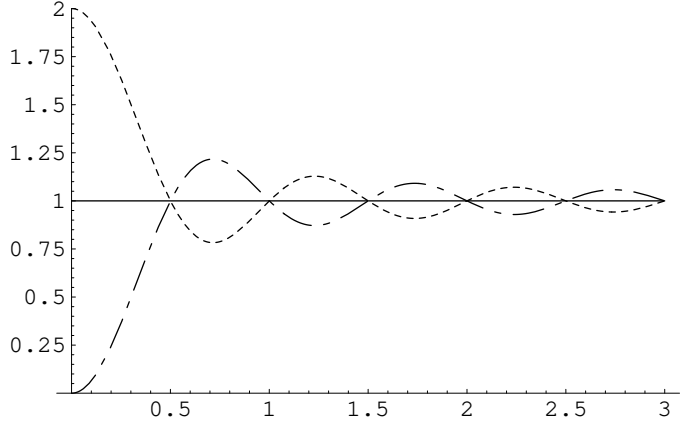}
  \caption{The central $1$-level densities $\rho_+$~(dotted) and
    $\rho_-$~(dash-dotted) versus the ``bulk'' $1$-level density
    $\rho\equiv1$ observed away from the central points.}
  \label{fig:1-level-U-O}
\end{figure}

\subsection{$1$-Level Density: Harder Orthogonal}
\label{sec:1-level-hard}

We return to the more general case of $m$-hard ensembles of orthogonal
matrices.  Because the classical Jacobi polynomials
$\{\pab_n\}_0^\infty$ are orthogonal with respect to the weight
$w(x)=(1-x)^a(1+x)^b$ on $[-1,1]$, the local statistics near the
central point $x=+1$ are derived from the asymptotic behavior of these
polynomials at the right edge of the interval $[-1,+1]$.%
\footnote{This ``edge limit'' and the ensuing Bessel Kernels are also
  observed in the somewhat simpler context of the so-called (unitary)
  Laguerre ensemble.}  More specifically, the relevant kernel which
takes the place of the (even or odd) Sine Kernel is the ``edge limit''
as $N\to\infty$ of the Christoffel-Darboux/Szeg\H{o} projection kernel
$K^{(a,b)}_N(x,y)$ onto polynomials of degree less than~$N$ in
$L^2([-1,1],(1-x)^a(1+x)^b dx)$ (via the change of
variables~\eqref{eq:xi}).  For the edge $+1$, the limit depends only
on the parameter~$a$ and is equal to the Bessel kernel%
\footnote{In fact, the even Sine Kernel $S_+=B^{(\mhalf)}$ whereas the
  odd Sine Kernel $S_-=B^{(\half)}$.}
\begin{align}
  \label{eq:2}
  \Ba(\xi,\eta) &= \frac{\sqrt{\xi\eta}}{\xi^2-\eta^2}
  [\pi\xi J_{a+1}(\pi\xi)J_a(\pi\eta) - J_a(\pi \xi) \pi\eta J_{a+1}(\pi \eta)], \\
  \label{eq:26}
  \Ba(\xi,\xi)
  &= \frac\pi2 (\pi\xi)[J_a^2(\pi\xi) - J_{a-1}(\pi\xi)J_{a+1}(\pi\xi)],
\end{align}
where $J_\nu$ stands for the Bessel function of the first kind~\cite{NW,D}.

It is a little more natural for our purposes to use the
hardness~$m$, rather than $a=m-\half$, as the parameter, so we
define
\begin{align}
  \label{eq:18}
  \Km(x,y) &= B^{(m\mhalf)}(x,y), \\
  \label{eq:25b}
  \km(x,y) &= \frac1\pi\Km(x/\pi,y/\pi).
\end{align}
Using the recursion formula for Bessel functions we obtain an
alternate formula to~\eqref{eq:26} for the diagonal values of the
kernel:
\begin{equation}
  \label{eq:17}
  \km(x,x) =
  \frac x2[J_{m+\half}(x)^2+J_{m-\half}(x)^2]
  - (m-{\textstyle\half})J_{m+\half}(x)J_{m-\half}(x).
\end{equation}
Except for $m$ times a point mass at $x=0$, the central $m$-hard
$1$-level density is given by
\begin{equation}
  \label{eq:31}
  \rho_m(x) = \Km(x,x).
\end{equation}
\begin{figure}[!htb]
  \centering
  \includegraphics[scale=.9]{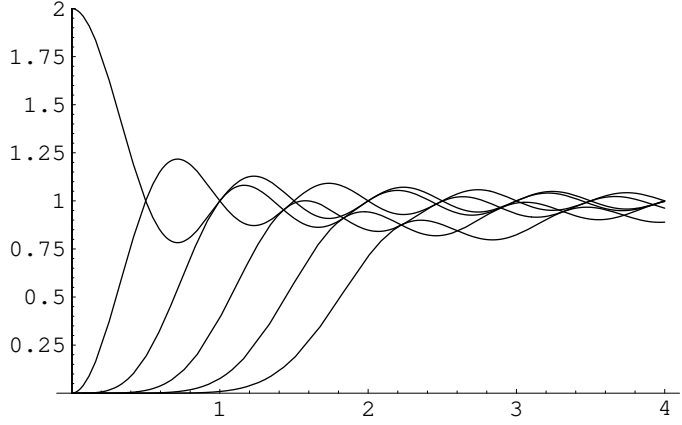}
  \caption{The $m$-hard $1$-level edge densities for $m=0,1,\dots,5$.}
  \label{fig:Bessel-1-density}
\end{figure}

\subsection{Spacing Measures}
\label{sec:spacing-meas} In this section we state some well-known
formulas giving the spacing measures or ``gap probabilities'' at
the central point.  Their derivation is standard and depends only
on knowledge of the edge limiting kernels~$\Km$ (see for instance
\cite{KaSa1,Meh,TW}). Let $E^{(m)}(k;s)$ be the limit, as
$N\to\infty$, of the probability that exactly $k$ of the $\xi_j$'s
lie on the interval $(0,s)$, where $\xi_j$ is related to $x_j$ via
equation~(\ref{eq:xi}). Also let $p^{(m)}(k;s)ds$ be the
conditional probability that the $(k+1)$-st of the $\xi_j$'s, to
the right of $\xi=0$, lies in the interval $[s,s+ds)$, in the
limit $N\to\infty$.

Abusing notation, let $\Km|_s$ denote the integral operator on
$L^2([0,s],dx)$ with kernel $\Km(x,y)$:
\begin{equation}
  \label{eq:32}
  \Km|_s f(\cdot) = \int_0^s \Km(\cdot,y)f(y)dy.
\end{equation}
If $I$ denotes the identity operator, then the following formulas
hold:
\begin{align}
  \label{eq:33}
  E^{(m)}(k;s) &= \frac1{k!}\left.\frac{\partial^k}{\partial T^k}
    \det(I+T\Km|_s)\right|_{T=-1}, \\
  \label{eq:34}
  p^{(m)}(k;s) &= -\frac d{ds}\sum_{j=0}^k E^{(m)}(j;s).
\end{align}
On the right-hand side of~(\ref{eq:33}), `det' is the Fredholm
determinant: for an operator with kernel $\mathcal{K}$,
\begin{equation}
  \label{eq:35}
  \det(I+\mathcal{K}) = 1 + \sum_{n=1}^\infty\frac1{n!}
  \idotsint_{\R^n}
  \limits\det_{n\times n}(\mathcal{K}(x_j,x_k)) \,dx_n\dots dx_1.
\end{equation}
Identical formulas hold even for finite $N$ provided that the
limiting kernel $K^{(m)}$ is replaced by the
Christoffel-Darboux/Szeg\H{o} projection kernel
$K^{(m-\half,\mhalf)}_{N-m}$ associated to the weight $w(x)$
of~(\ref{eq:11}) with $a=m-\half$, acting on $L^2([-1,1],w(x)dx)$.
In this case, the corresponding operator is of finite rank, the
Fredholm determinant agrees with the usual determinant, and the
series~(\ref{eq:35}) is finite.

\subsection{Explicit Kernels}
\label{sec:explicit-kernels} In view of the relation between
Bessel Functions of the first kind of half-integral parameter and
trigonometric functions, it is possible to write the kernels $\Km$
in terms of elementary functions.
\subsubsection{$m=0$: The Even Sine Kernel}
\label{sec:even-sine}
\begin{equation}
  \label{eq:19}
  K_{0}(x,y) = S_+(x,y)
  = \frac{\sin\pi(x-y)}{\pi(x-y)} + \frac{\sin\pi(x+y)}{\pi(x+y)}.
\end{equation}
The one-level density is
\begin{equation}
  \label{eq:20}
  \rho_+(x) = S_+(x,x) = 1 + \frac{\sin2\pi x}{2\pi x}.
\end{equation}
The Fourier transform of the one-level density is
\begin{equation}
  \label{eq:21}
  \hat\rho_+(u) = \delta(u) + \half I(u),
\end{equation}
where $I(u)$ is the characteristic function of the interval
$[-1,1]$.

\subsubsection{$m=1$: The Odd Sine Kernel}
\label{sec:odd-sine}
\begin{equation}
  \label{eq:22}
  K_{1}(x,y) = S_-(x,y)
  = \frac{\sin\pi(x-y)}{\pi(x-y)} - \frac{\sin\pi(x+y)}{\pi(x+y)}.
\end{equation}
The one-level density is
\begin{equation}
  \label{eq:23b}
  \rho_-(x) = S_-(x,x) = \delta(x) + 1 - \frac{\sin2\pi x}{2\pi x}.
\end{equation}
The Fourier transform of the one-level density is
\begin{equation}
  \label{eq:24b}
  \hat\rho_-(u) = \delta(u) + 1 - \half I(u).
\end{equation}

\subsubsection{$m=2$: The ``Doubly Hard'' Kernel}
\label{sec:2-hard}
\begin{equation}
  \label{eq:25c}
  K_{2}(x,y)
  = \frac{\sin\pi(x-y)}{\pi(x-y)} + \frac{\sin\pi(x+y)}{\pi(x+y)}
  - 2\frac{\sin \pi x}{\pi x}\frac{\sin \pi y}{\pi y}.
\end{equation}
The one-level density is
\begin{equation}
  \label{eq:23c}
  \rho_2(x)
  = 2\delta(x) + 1 + \frac{\sin2\pi x}{2\pi x} - 2\left(\frac{\sin\pi x}{\pi x}\right)^2.
\end{equation}
The Fourier transform of the one-level density is
\begin{equation}
  \label{eq:24c}
  \hat\rho_2(u) = \delta(u) + 2 + (2|u|-{\textstyle\frac32}) I(u).
\end{equation}

\subsubsection{$m=3$: The ``Triply Hard'' Kernel}
\label{sec:3-hard}
\begin{equation}
  \label{eq:27}
  \begin{split}
    K_3(x,y) = K_1(x,y) +
    \frac{18}{\pi^2xy}\left(1+\frac5{\pi^2xy}\right)K_0(x,y)\\
    {}- 6\left(\frac{\cos\pi x}{\pi x}\frac{\cos\pi y}{\pi y}
      +\frac{\sin\pi x}{(\pi x)^2}\frac{\sin\pi y}{(\pi y)^2}\right).
  \end{split}
\end{equation}

\end{document}